\documentclass[11pt,english,a4paper]{smfart}
\usepackage[english]{babel}
\usepackage{xcolor}
\usepackage{pbsi}
\usepackage[T1]{fontenc}
\usepackage{stmaryrd}
\usepackage{mathabx}
\usepackage[mathscr]{euscript}
\usepackage{amsmath}
\usepackage[pagebackref,colorlinks=true,linkcolor=blue,citecolor=blue,urlcolor=red, hypertexnames=true]{hyperref}
 \usepackage[abbrev,backrefs,nobysame]{amsrefs}
 \usepackage{caption}
 \usepackage{blkarray}
 \usepackage{enumerate}
\usepackage{amsmath,amssymb,bm,amsfonts}
\usepackage[all]{xy}
\entrymodifiers={+!!<0pt,\fontdimen22\textfont2>}

\usepackage{mathrsfs}

\usepackage{graphicx}

\usepackage{color}

\makeatletter
\renewcommand\theequation%
{\thesection.\arabic{equation}}
\@addtoreset{equation}{section}
\makeatother

\newtheorem{theorem}{Theorem}[section]
\newtheorem{lemma}[theorem]{Lemma}

\newtheorem{corollary}[theorem]{Corollary}
\newtheorem{proposition}[theorem]{Proposition}

\newtheorem*{coro*}{Corollary}

{\theoremstyle{definition}}

{\theoremstyle{definition}\newtheorem{example}[theorem]{Example}}

\theoremstyle{definition}
\newtheorem{definition}[theorem]{Definition}
\newtheorem{question}[theorem]{Question}
\newtheorem{fact}[theorem]{Fact}

\newtheorem*{remark*}{Remark}

{\theoremstyle{definition}\newtheorem{remark}[theorem]{Remark}}

\def\K{\ensuremath{\mathbb K}}

\def\R{\ensuremath{\mathbb R}}
\def\Z{\ensuremath{\mathbb Z}}

\def\C{\ensuremath{\mathbb C}}

\def\N{\ensuremath{\mathbb N}}

\def\bth{\begin{theorem}}
\def\blm{\begin{lemma}}
\def\bpr{\begin{proposition}}
\def\bpf{\begin{proof}}
\def\epf{\end{proof}}
\def\epr{\end{proposition}}
\def\elm{\end{lemma}}
\def\eth{\end{theorem}}
\def\bco{\begin{corollary}}
\def\eco{\end{corollary}}
\def\be{\begin{enumerate}}
\def\ee{\end{enumerate}}
\def\bea{\begin{enumerate}[\rm (a)]}
\def\beun{\begin{enumerate}[\rm (1)]}
\def\bei{\begin{enumerate}[\rm (i)]}

\def\bdf{\begin{definition}}
\def\edf{\end{definition}}

\def\u{{\pmb u}}
\def\v{{\pmb v}}
\def\w{{\pmb w}}

\def\vprod{{v_1\cdots v_n}}
\def\uprod{{u_1\cdots u_n}}
\def\wprod{{w_1\cdots w_n}}

 \makeatletter
\@namedef{subjclassname@2020}{%
  \textup{2020} Mathematics Subject Classification}
\makeatother

\author[S. Grivaux]{Sophie Grivaux}
\address[S. Grivaux]{CNRS, Univ. Lille, UMR 8524 - Laboratoire Paul
Painlev\'e, F-59000 Lille, France}
\email{sophie.grivaux@univ-lille.fr}
\author[\'{E}. Matheron]{\'{E}tienne Matheron}
\address[\'{E}. Matheron]{Univ. Artois, UR 2462 - Laboratoire de Math\'{e}matiques de Lens (LML)\\ F-62300 Lens, France}
\email{etienne.matheron@univ-artois.fr}
\author[Q. Menet]{Quentin Menet}
\address[Q. Menet]{Service de Probabilit\'e et Statistique, D\'epartement de Math\'ematique\\ Universit\'{e} de Mons\\ Place du Parc 20\\ 7000 Mons, Belgium}
\email{quentin.menet@umons.ac.be}

\title[Invariant measures for weighted shifts]{Orthogonality of invariant measures\\ for weighted shifts}

\begin{abstract}
We introduce and study the notion of \emph{orthogonality} for two operators in the context of weighted backward shifts on $\ell_p(\Z_+)$, $1\leq p<\infty$. Two continuous linear operators $T_1$ and $T_2$ acting on a Polish topological vector space $X$ are said to be {orthogonal} if any two Borel probability measures $m_1$ and $m_2$ on $X$ which are respectively $T_1\,$-$\,$invariant and $T_2\,$-$\,$invariant and  satisfy $m_1(\{0\})=m_2(\{0\})=0$ must be orthogonal. In this note, we provide several conditions on the weights $\u$ and $\v$ implying orthogonality or non-orthogonality of the associated weighted shifts $B_\u$ and $B_\v$, and we investigate in some detail the case where the invariant measures are product measures.
\end{abstract}

 \thanks{The third author is a Research Associate of the Fonds de la Recherche Scientifique - FNRS}

\begin{document}
\maketitle

\hfill\emph{To Gilles Godefroy, with affection and admiration}

\bigskip

\section{Introduction}\label{Intro}
\subsection{Weighted shifts} This note is a contribution to the study of the dynamics of weighted backward shifts acting on $\ell_p\,$-$\,$spaces, $1\leq p<\infty$. More specifically, we will be concerned with \emph{invariant measures} for such weighted shifts.

\smallskip
Let $X$ be one of the (real or complex) Banach spaces $\ell_p=\ell_p(\Z_+)$, $1\leq p<\infty$, and denote by $(e_n)_{n\geq 0}$ the canonical basis of $X$. Let $\pmb w=(w_n)_{n\geq 1}$ be a \emph{weight sequence}, \textit{i.e.} a bounded sequence of non-zero scalars.The backward shift associated with $\pmb w$ is the operator $B_{\pmb w}:X\to X$ defined by $B_{\pmb w}e_0:=0$ and $B_{\pmb w}e_n:=w_n e_{n-1}$ for $n\geq 1$. Weighted shift operators occupy a central place in operator theory. On the one hand, their explicit form permits a systematic study, so that they are natural ``test'' operators for any question that comes to mind. On the other hand, their theory is quite rich, and their behaviour remains in some aspects rather mysterious. We refer the reader to the classical paper \cite{Sh}, which proposes a systematic study of weighted shifts including their representations as multiplication operators on certain Banach spaces of holomorphic functions and their spectral properties, and to the authoritative book \cite{Ni}, which deals with the study of cyclic vectors for the unweighted backward shift $B$ on $\ell_2(\Z_+)$. 

\smallskip
Many important dynamical properties are completely characterized for weighted shifts $B_\w$ and can be expressed in a simple way in terms of the weights $w_n$, see e.g \cite{Sal1}, \cite{Sal2}, \cite{CS}, \cite{BCDMP}, \cite{BerMes}.  For instance,  $B_{\pmb w}$ is \emph{hypercyclic}, \textit{i.e.} it admits a vector with a dense orbit, if and only if 
$\overline\lim_{n\to\infty}\, \vert w_1\cdots w_n\vert=\infty$, and it is \emph{topologically mixing} if and only if $\vert\wprod\vert\to\infty$ as $n\to\infty$. 
See \cite{GEP} for proofs of these results as well as extensions to a large class of sequence spaces. Observe that these characterizations are the same on any $\ell_p$. By contrast,  \emph{chaoticity} of a weighted shift on $X=\ell_p$ depends explicitely on $p$ (recall that an operator is chaotic if it is hypercyclic with a dense set of periodic points): $B_{\pmb w}$ is chaotic on $\ell_p$ if and only $\sum_{n=1}^\infty \frac1{\vert w_1\cdots w_n\vert^p}<\infty$. For the class of weighted shift operators of $\ell_p$, chaos is an especially important notion since it turns out to be equivalent to strong dynamical properties such as \emph{frequent hypercyclicity}  (the existence of a vector whose orbit visits each non-empty open set along a set of integers having positive lower  density) or \emph{$\mathcal U$-$\,$frequent hypercyclicity}  (same as frequent hypercyclicity with lower density replaced by upper density). Indeed, by an important result of Bayart and Ruzsa \cite{BR}, frequent hypercyclicity and $\mathcal U$-$\,$frequent hypercyclicity of a weighted shift $B_{\pmb w}$ on $\ell_p$ are both equivalent to the condition $\sum_{n=1}^\infty \frac1{\vert w_1\cdots w_n\vert^p}<\infty$. 
 We refer the reader to the books \cite{BM} and \cite{GEP} for more on frequent hypercyclicity and related questions, and to the paper \cite{ChGEMe} for a study of the relation between chaos and frequent hypercyclicity for weighted shifts on a large class of Fr\' echet sequence spaces.
 
 \smallskip
 The condition $\sum_{n=1}^\infty \frac1{\vert w_1\cdots w_n\vert^p}<\infty$ also turns out to be of special importance in the study of {invariant measures} for the weighted backward shift $B_{\pmb w}$. Recall that a Borel probability measure $m$ on $X$ is said to be $B_{\pmb w}\,$-$\,$invariant if 
 $m\bigl( B_{\pmb w}^{-1}(A)\bigr)=m(A)$ for every Borel set $A\subseteq X$. Under the condition $\sum_{n=1}^\infty \frac1{\vert w_1\cdots w_n\vert^p}<\infty$, the operator $B_{\pmb w}$ acting on $X=\ell_p$ admits plenty of invariant measures, of quite different kinds: discrete measures associated with periodic points, Gaussian measures with full support (with respect to some of which $B_{\pmb w}$ is ergodic or even strongly mixing), continuous measures with full support which are very far from being Gaussian, ... See \cite{BG1} or \cite{BM}*{Chapter 5} for more information on (shift-invariant) Gaussian measures; and \cite{G}, which testifies of the richness of the class of invariant measures for backward shifts by showing that under the condition $\sum_{n=1}^\infty \frac1{\vert w_1\cdots w_n\vert^p}<\infty$, the backward shift $B_{\pmb w}$ is \emph{universal for ergodic system} in the sense of Glasner and Weiss \cite{GW}: for every ergodic transformation $T$ on a standard Lebesgue probability space $(Z,\mathcal B, \mu)$, there exists a $B_{\pmb w}\,$-$\,$invariant Borel probability measure $m$ on $X=\ell_p$ with full support such that the two dynamical systems $(Z,\mathcal B,\mu ;T)$ and $(X,\mathcal B_X, m;B_{\pmb w})$ are isomorphic (where $\mathcal B_X$ is the Borel $\sigma$-algebra of $X$).
 
 \smallskip Invariant measures for weighted backward shifts are thus far from being completely understood (and in some sense, they will never be). In this note, our aim is to contribute to their understanding by studying possible links between invariant measures for \emph{two different} weighted shifts $B_{\pmb u}$ and $B_{\pmb v}$. More specifically, we will be concerned with the notion of \emph{orthogonality} of two weighted shifts.
 
 \smallskip 
 \subsection{Orthogonality} Before giving the definition we will be playing with, we fix some notation. The scalar field is denoted by $\K$. All measures on $X=\ell_p(\Z_+)$ are understood to be Borel probability measures. Given a weight sequence $\pmb w$, we denote by $\mathcal P_{\pmb w}(X)$ the set of all $B_{\pmb w}\,$-$\,$invariant measures on $X$. If $m$ and $m'$ are two measures on $X$, we write $m\perp m'$ to indicate that $m$ and $m'$ are orthogonal (\textit{i.e.} mutually singular), and $m\ll m'$ to indicate that $m$ is absolutely continuous with respect to $m'$. If $m\ll m'$ and $m'\ll m$, we say that $m$ and $m'$ are equivalent and we write $m\sim m'$. 
 
 \smallskip We now define orthogonality of two weighted backward shifts.
 \bdf\label{defortho} Two weighted backward shifts $B_{\pmb u}$ and $B_{\pmb v}$ on $X=\ell_p$ are said to be \emph{orthogonal} if the following holds true: whenever $m_{\pmb u}\in \mathcal P_{\pmb u}(X)$, $m_{\pmb v}\in \mathcal P_{\pmb v}(X)$ and $m_{\pmb u}(\{ 0\})=0=m_{\pmb v}(\{ 0\}) $, it follows that $m_{\pmb u}\perp m_{\pmb v}$. 
  \edf
 
\smallskip Removing the $\delta_0\,$-$\,$parts of the measures, we see that this definition can be reformulated as follows: $B_\u$ and $B_\v$ are orthogonal if and only if whenever $m_{\pmb u}\in \mathcal P_{\pmb u}(X)$ and $m_{\pmb v}\in \mathcal P_{\pmb v}(X)$ are non-orthogonal, it must be  that $m_\u(\{ 0\})>0$ \emph{and} $m_\v(\{0\})>0$. Equivalently: there exists no (Borel probability) measure on $X$ which is absolutely continuous with respect to both a $B_\u\,$-$\,$invariant measure and a $B_\v\,$-$\,$invariant measure, except the Dirac mass $\delta_0$. Thus, 
informally speaking, $B_{\pmb u}$ and $B_{\pmb v}$ are orthogonal if their invariant measures ``have nothing to say to each other'', unless they charge the singleton $\{ 0\}$. The main question we want to investigate is the following.
 \begin{question}\label{Q1} Characterize the pairs of weight sequences $(\pmb u,\pmb v)$ such that the associated weighted shifts $B_{\pmb u}$ and $B_{\pmb v}$ acting on $X=\ell_p$ are orthogonal.
 \end{question} 
 
 \smallskip The motivation for looking at this question comes from a recent work of Charpentier, Ernst, Mestiri and Mouze \cite{CEMM}. Let $B$ be the unweighted backward shift acting on the complex Hilbert space $\ell_2(\Z_+)$, and let $\Lambda\subseteq\C$. In \cite{CEMM}, it is shown that the operators $\lambda B$, $\lambda\in\Lambda$ admit a \emph{common frequently hypercyclic vector} if and only if the set $\{ \vert \lambda\vert;\; \lambda\in\Lambda\}$ is a countable relatively compact subset of $(1,\infty)$. In particular, if $a,b\in\C$ and $\vert a\vert, \vert b\vert >1$, then the operators $aB$ and $bB$ have a common frequently hypercyclic vector. This is by no means obvious if $\vert a\vert\neq \vert b\vert$. So it is tempting to ask whether this result could be retrieved in a ``soft'' way by measure-theoretic arguments. Indeed, imagine that it were possible to find two measures $m_a$ and $m_b$ on $\ell_2$, with full support, invariant and \emph{ergodic} for $aB$ and $bB$ respectively, such that one of them is absolutely continuous with respect to the other, say $m_a\ll m_b$. By the pointwise  ergodic theorem, it would follow that $m_a\,$-$\,$almost every $x\in\ell_2$ is frequently hypercyclic for $aB$ and $m_b\,$-$\,$almost every $x\in\ell_2$ is frequently hypercyclic for $bB$; hence $m_a\,$-$\,$almost every $x$ would be frequently hypercyclic for both $aB$ and $bB$ since $m_a\ll m_b$, so one could conclude in particular that $aB$ and $bB$ have a common frequently hypercyclic vector. However, such measures $m_a$ and $m_b$ simply do not exist: we will see below that the operators $aB$ and $bB$ are in fact orthogonal.
 
 \subsection{Organization of the paper} In Section \ref{Sec2}, we complement existing results by showing that a weighted backward shift $B_\w$ on $X=\ell_p(\Z_+)$ admits a non-trivial invariant measure \emph{if and only if} $\sum_{n=1}^\infty \frac1{\vert\wprod\vert^p}<\infty$. In Section \ref{Sec3}, we define orthogonality in a very general context, we give a simple condition ensuring orthogonality (Theorem \ref{Th1}), and we use this to show, among other things, that if $T$ is any continuous linear operator acting on a Polish topological vector space $X$ and $a, b\in\K$ are such that $\vert a\vert\neq\vert b\vert$, then the operators $aT$ and $bT$ are orthogonal (Corollary \ref{ab}). We also deduce from our criterion that two weighted shifts that  are ``far from being similar'' must be orthogonal (Theorem \ref{presquebis}). In Section \ref{Sec4}, we use \emph{periodic vectors} to give examples of non-orthogonal weighted shifts; in particular, we show that there exist non-orthogonal weighted shifts which are not similar (Example \ref{non}). Section \ref{Sec5} is devoted to invariant \emph{product measures} for weighted shifts acting on $\ell_p$, and their role in our study of orthogonality. We give a necessary and sufficient condition for two weighted shifts $B_\u$ and $B_\v$ to admit equivalent invariant Gaussian product measures (Corollary \ref{gauss'}), and then we examine to what extent the existence of non-orthogonal invariant product measures implies the existence of equivalent Gaussian product measures (Theorem \ref{gauss2} and Proposition \ref{gauss3}). Section \ref{add} contains some  additional facts. Some of these facts are actually used in earlier proofs, but we prefered to postpone them for the sake of fluid reading. We conclude the paper with a few natural questions. 
 
 \section{Existence of non-trivial invariant measures}\label{Sec2}
 It is, of course, not very interesting to study orthogonality of two weighted shifts if one does not know at least that each of them admits non-trivial invariant measures, \textit{i.e.} invariant measures different from the Dirac mass $\delta_0$. The following proposition says exactly when this happens for a weighted shift acting on $\ell_p$.
 \bpr\label{existe} Let $B_\w$ be a weighted backward shift acting on $X=\ell_p$. Then $B_\w$ admits a non-trivial invariant measure if and only if $\sum\limits_{n=1}^\infty \frac1{\vert w_1\cdots w_n\vert^p}<\infty$.
 \epr
 
 One can give a short proof of Proposition \ref{existe} by using a (quite non-trivial) lemma inspired by \cite{BR},  which can be extracted from the proof of \cite{ChGEMe}*[Theorem 2.1]. Since this lemma might be useful in other situations, we state it explicitely. Let us denote by $e_0^*$ the first coordinate functional on $X=\ell_p$. 
 
 \blm\label{summability} Let $B_\w$ be a weighted shift acting on $X=\ell_p(\Z_+)$. Assume that one can find a vector $x\in\ell_p$ and a set of integers $\mathcal N$ such that
 \be
 \item[\rm (i)] $\mathcal N$ has positive upper density,
 \item[\rm (ii)] $\sup_{n\in\mathcal N} \Vert B_\w^n x\Vert<\infty$,
 \item[\rm (iii)] $\inf_{n\in\mathcal N} \vert \langle e_0^*, B_\w^n x\rangle\vert >0$.
 \ee
 
 Then, one can conclude that $\sum\limits_{n=1}^\infty \frac1{\vert w_1\cdots w_n\vert^p}<\infty$.
 \elm
  \begin{proof}[Sketch of the proof of Lemma \ref{summability}] Write $x=\sum_{n=0}^\infty x_n e_n$. 
Let $C_1:=\sup_{n\in\mathcal N} \Vert B_\w^n x\Vert$ and $C_2:=\inf_{n\in\mathcal N} \vert \langle e_0^*, B_\w^n x\rangle\vert=\inf_{n\in\mathcal N} \vert \wprod\, x_n\vert$.
If we let $v_{n}:=(w_1\cdots w_n)^{-1}$, we have for every $n\in \mathcal{N}$ and every integer $M\ge 0$,
\begin{align*}
\Biggl\Vert\sum_{\substack{m\in \mathcal{N}\\n\leq m\leq n+M}}v_{m-n}e_{m-n}\Biggr\Vert &\leq \frac{1}{C_2}\Biggl\Vert\sum_{\substack{m\in \mathcal N\\n\leq m\leq n+M}}v_{m-n}\frac{x_m}{v_m}e_{m-n}\Biggr\Vert= \frac{1}{C_2}\Biggl\Vert\sum_{\substack{0\leq k\leq M\\k+n\in \mathcal N}}v_{k}\frac{x_{k+n}}{v_{k+n}}e_{k}\Biggr\Vert\leq \frac{C_1}{C_2}\cdot
\end{align*}

We then conclude thanks to \cite{ChGEMe}*[Lemma 2.5].
\end{proof}
 
 \smallskip
 \begin{proof}[Proof of Proposition \ref{existe}]  As mentioned in the introduction, it is well-known that if $\sum_{n=1}^\infty \frac1{\vert w_1\cdots w_n\vert^p}<\infty$, then $B_\w$ admits lots of non-trivial invariant measures. Perhaps the simplest such measure is the Dirac mass $\delta_x$ at the fixed point $x:=e_0+\sum_{n=1}^\infty \frac1\wprod\, e_n$.
 
 \smallskip Conversely, assume that  $B_\w$ admits a non-trivial invariant measure $m$. By the ergodic decomposition theorem, we may assume that $m$ is an ergodic measure for $B_\w$. Since $m\neq \delta_0$, we may choose a point $u\neq 0$ in the support of $m$. Then $B_\w^nu\in{\rm supp}(m)$ for all $n\geq 0$ by the $B_\w\,$-$\,$invariance of $m$, so we may assume that $\langle e_0^*,u\rangle\neq 0$. Let $\alpha>0$ be such that $\vert \langle e_0^*,u\rangle\vert >\alpha$, and consider the open set 
 \[ U:=\bigl\{ z\in X;\; \Vert z-u\Vert<1\quad\hbox{and}\quad\vert \langle e_0^*,z\rangle\vert>\alpha\bigr\}.\]
 Since $U$ is a neighbourhood of $u$, we have $m(U)>0$. By the pointwise ergodic theorem, it follows that one can find $x\in X$ such that the set $\mathcal N_{B_\w}(x,U):=\{ n\in\Z_+;\; B_\w^nx\in U\}$ admits a positive density, \textit{i.e.} 
 $
 \lim_{N\to\infty} \frac1{N+1}\# \bigl( \llbracket 0,N\rrbracket\cap \mathcal N_{B_\w}(x,U)\bigr)$ exists and is (strictly) positive; in particular, $\mathcal N_{B_\w}(x,U)$ has positive \emph{upper} density.  So we have found a vector $x\in\ell_p$ and a set of integers $\mathcal N$ with positive upper density such that $\Vert B_\w^n x\Vert<1+\Vert u\Vert$ and $\vert \langle e_0^*, B_\w^n x\rangle\vert>\alpha$ for all $n\in\mathcal N$.
  By Lemma \ref{summability}, this concludes the proof. 
 \epf
 
 \begin{remark} Proposition \ref{existe} is valid in a more general setting: the same proof shows that if $X$ is a Banach sequence space for which $(e_n)_{n\geq 0}$ is a  {boundedly complete} unconditional basis, then a weighted shift $B_\w$ acting on $X$ admits a non-trivial invariant measure if and only if the series $\sum \frac1\wprod\, e_n$ is convergent in $X$. {We do not know if the analogous result holds for weighted shifts acting on $c_0(\Z_+)$.}
 \end{remark}
 
 
\section{A general criterion for orthogonality}\label{Sec3} Although the definition of orthogonality was given in the introduction for pairs of weighted shifts only, it makes sense in a much more general context. If $T$ is any continuous self-map of a Polish (\textit{i.e.} separable and completely metrizable) space $X$, let us denote by $\mathcal P_T(X)$ the set of all $T\,$-$\,$invariant (Borel probability) measures on $X$.
\bdf Let $X$ be a Polish space, and let {\rm (P)} be any property of measures on $X$. We say that two continuous self-maps $T_1, T_2$ of $X$  are \emph{orthogonal with respect to measures satisfying {\rm (P)}} if the following holds true: whenever $m_1\in\mathcal P_{T_1}(X)$ 
and $m_2\in\mathcal P_{T_2}(X)$ satisfy {\rm (P)}, it follows that $m_1\perp m_2$.
\edf

In accordance with  Definition \ref{defortho}, if $X$ is a Polish topological vector space and $T_1,T_2$ are continuous linear operators on $X$, we will say that $T_1$ and $T_2$ are \emph{orthogonal} if they are orthogonal with respect to measures not charging the singleton $\{ 0\}$.

\smallskip
We now present a general sufficient condition for orthogonality. Recall that if $(A_n)_{n\in\N}$ is a sequence of subsets of $X$, then $\varlimsup\, A_n$ denotes the set of all $x\in X$ that belong to infinitely many $A_n$'s.
\bth\label{Th1} Let $X$ be a Polish space, and let $\Lambda$ be a Borel subset of $X$. Let also $T_1$ and $T_2$ be two continuous self-maps of $X$. Assume that for any compact set $K\subseteq X$ with 
$K\cap \Lambda=\emptyset$, there exists an infinite set $I\subseteq \N$ such that 
\begin{equation}\label{limsup0} \overline{\lim_{n\in I}} \;\, T_1^{-n}(K)\cap T_2^{-n}(K)=\emptyset .
\end{equation}

Then, $T_1$ and $T_2$ are orthogonal with respect to measures not charging $\Lambda$.
\eth
\bpf Let $m_1\in\mathcal P_{T_1}(X)$ and $m_2\in\mathcal P_{T_2}(X)$ be such that $m_1(\Lambda)=0=m_2(\Lambda)$. We have to show that $m_1\perp m_2$, and this will be an easy consequence of the following fact.
\begin{fact}\label{factCharge} For any $\varepsilon >0$, one can find a Borel set $E_\varepsilon\subseteq X$ such that $m_1(E_\varepsilon)<\varepsilon$ and $m_2(E_\varepsilon)>1-\varepsilon$.
\end{fact}
\begin{proof}[Proof of Fact \ref{factCharge}] By regularity of the measures $m_1$ and $m_2$ and since $m_1$ and $m_2$ do not charge $\Lambda$, there exists a compact set $K\subseteq X\setminus\Lambda$ such that $m_1(K)>1-\varepsilon/2$ and $m_2(K)>1-\varepsilon/2$. Let $I$ be an infinite subset of $\N$ such that (\ref{limsup0}) holds true. Then $m_1\bigl( T_1^{-n}(K)\cap T_2^{-n}(K)\bigr)\to 0$ as $n\to\infty$ along $I$ by Fatou's Lemma; so one can find an integer $n\geq 1$ such that $m_1\bigl( T_1^{-n}(K)\cap T_2^{-n}(K)\bigr)<\varepsilon/2$. Since $m_1\bigl( T_1^{-n}(K)\bigr)=m_1(K)>1-\varepsilon/2$, it follows that 
$m_1\bigl( T_2^{-n}(K)\bigr)<\varepsilon$. On the other hand, $m_2\bigl( T_2^{-n}(K)\bigr)=m_2(K)>1-\varepsilon$; so we may take $E:= T_2^{-n}(K)$.
\epf

If we now define $E:=\bigcup_{k\in\N} \bigcap_{n\in \N} E_{{\frac{2^{-n}}k}}$, we see that  
 $m_1(E)=0$ and $m_2(E)=1$, which shows that $m_1\perp m_2$.
\epf

\begin{remark} The assumption in Theorem \ref{Th1} is here to ensure that for any compact set $K\subseteq X\setminus\Lambda$ and any measure $m$ on $X$, it holds that $\inf_{n\geq 1} m\bigl( T_1^{-n}(K)\cap T_2^{-n}(K)\bigr)=0$. However, it is very likely that (\ref{limsup0}) is a too strong assumption. More generally, it would be nice to have a ``measure-free'' characterization of the sequences of Borel sets $(A_n)_{n\geq 1}\subseteq X$ such that $\inf_{n\geq 1} m(A_n)=0$ for every measure $m$ on $X$. Note that there is a very simple characterization if ``$\inf$'' is replaced by ``$\lim$'': we have $\lim_{n\to\infty} m(A_n)=0$ for every measure $m$ on $X$ if and only if $\varlimsup A_n=\emptyset$. On the other hand, for a sequence $(A_n)_{n\geq 1}$ of non-empty Borel subsets of $X$, consider the following statements.
\be
\item[\rm (i)] There exists an infinite set $I\subseteq \N$ such that $\varlimsup_{n\in I} A_n=\emptyset$.
\item[\rm (ii)] $\inf_{n\geq 1} m(A_n)=0$ for every measure $m$ on $X$.
\item[\rm (iii)] $\inf_{n\geq 1} m(A_n)=0$ for every discrete measure $m$ on $X$.
\item[\rm (iv)] $\inf_{n\geq 1} m(A_n)=0$ for every finitely supported measure $m$ on $X$.
\item[\rm (v)] For any finite set $F\subseteq X$, there exist infinitely many $n$ such that $A_n\cap F=\emptyset$.
\ee

Then (i)$\implies$(ii)$\implies$(iii)$\iff$(iv)$\iff$(v), but (v) does not imply (ii) and (ii) does not imply (i).
\end{remark}
\bpf (i)$\implies$(ii) by Fatou's Lemma, and obviously (ii)$\implies$(iii)$\implies$(iv). Moreover, it is rather clear that in fact (iii)$\iff$(iv), and that (v)$\implies$(iv).

\smallskip
To prove that (iv)$\implies$(v), assume that for some finite set $F\subseteq X$ and some integer $N$, we have $A_n\cap F\neq\emptyset$ for all $n>N$. Choose a point $a_i\in A_i$ for $i=1,\dots ,N$. Then, the finitely supported measure $m:=\frac12\left(\frac1{N}\sum\limits_{i=1}^N \delta_{a_i}+\frac1{\#F}\sum\limits_{x\in F} \delta_x\right)$ is such that $\inf_{n\geq 1} m(A_n)>0$.


\smallskip To show that (v) does not imply (ii), let $X:= [0,1]$, and define \[ A_n:=\{ t\in [0,1];\; \vert \sin(2\pi n t)\vert \geq 1/\pi\}.\]

\noindent By Dirichlet's theorem, for any finite set $F\subseteq [0,1]$ one can find an increasing sequence of integers $(n_k)$ such that $\sin(2\pi n_k t)\to 0$ on $F$; so (v) is satisfied. However, we have $\int_0^1 \vert \sin(2\pi nt)\vert \, dt=2/\pi$ for all $n\geq 1$. Hence $m(A_n)\geq 1/\pi$ for all $n\geq 1$, where 
$m$ is  Lebesgue measure on $[0,1]$.

\smallskip Let us now show that (ii) does not imply (i). This example is due to N. de Rancourt. Consider the Cantor space $\Delta=\{ 0,1\}^\N$ identified with the power set of $\N$. Let 
\[ X:=\left\{ J\subseteq\N;\; \sum_{n\in J} \frac1n\leq 1\right\},\]
which is a closed subset of $\Delta$ and hence a (compact) Polish space. For each $n\in\N$, let 
\[ A_n:=\{ J\in X;\; n\in J\}.\]

It is clear that the sequence $(A_n)$ does not satisfy (i): indeed, for any infinite set $I\subseteq\N$, one can find an infinite set $J\subseteq I$ such that $\sum_{n\in J}\frac1n\leq 1$, and this $J$ belongs to $\varlimsup A_n$ by definition.

However, $(A_n)$ satisfies (ii); in fact, for any measure $m$ on $X$, we have $\sum_{n=1}^\infty \frac{m(A_n)}n<\infty$. Indeed, we have 
\[ \sum_{n=1}^\infty \frac{m(A_n)}n=\int_X\left( \sum_{n=1}^\infty \frac1n \,\mathbf 1_{A_n}(J)\right) dm(J)=\int_X\left(\sum_{n\in J} \frac1n\right) dm(J)\leq 1.\]


\smallskip Here is yet another example showing that (ii) does not imply (i), which is a kind of twin of the previous one. Let
\[ X:=\left\{ x\in [0,\pi];\; \sum_{n=1}^\infty \frac{\vert \sin(2^nx)\vert}{n}\leq 1\right\},\]
which is a closed subset of $[0,\pi]$. For each $n\geq 1$, let
\[ A_n:=\{ x\in X;\; \vert \sin(2^n x)\vert \geq 1/2\}.\]

With some effort, one can prove the following: \emph{if $(m_i)_{i\geq 1}$ is a sufficiently fast increasing sequence of integers and if we define
\[ x:=\frac{\pi}{2} \,\sum_{i=1}^\infty  \frac1{2^{m_i}},\]
then $x\in A_{m_k}$ for every $k\geq 1$.} It follows imediately that the sequence $(A_n)$ does not satisfy (i). However, one shows in the same way as above that $(A_n)$ satisfies (ii).
\epf


\begin{remark}\label{unpeuutile} With the notation of Theorem \ref{Th1}, consider the following statements.
\be
\item[(1)] For any compact set $K\subseteq X\setminus\Lambda$, it holds that $\overline{\lim} \; T_1^{-n}(K)\cap T_2^{-n}(K)=\emptyset$.
\item[(2)] It is not possible to find $x\in X$ and an increasing sequence of integers $(n_k)$ such that the sequences $(T_1^{n_k}x)$ and $(T_2^{n_k}x)$ both converge to a limit not belonging to $\Lambda$.
\ee

Then (2)$\implies$(1), and (1)$\iff$(2) if $\Lambda$ is a \emph{closed} subset of $X$.
\end{remark}
\bpf It should be rather clear that (2) implies (1). Indeed, let $K$ be a compact subset of $X\setminus \Lambda$, and assume that $\overline{\lim} \; T_1^{-n}(K)\cap T_2^{-n}(K)\neq\emptyset$. Then, there exists $x\in X$ and an increasing sequence of integers $(n_k)_{k\geq 0}$ such that $T_1^{n_k}x\in K$ and 
$T_ 2^{n_k}x\in K$ for all $k\geq 0$. Since $K$ is compact, we may assume, upon extracting subsequences, that both sequences $(T_1^{n_k}x)$ and $(T_2^{n_k}x)$ converge, $T_1^{n_k}x\to u\in K$ and $T_2^{n_k}x\to v\in K$; and since $K\subseteq X\setminus\Lambda$, this shows that (2) is not satisfied.

Conversely, assume that $\Lambda$ is closed and that (2) is not satisfied, \textit{i.e.} there exist $x\in X$ and an increasing sequence of integers $(n_k)_{k\geq 0}$ such that $T_1^{n_k}x\to u$ and $T_2^{n_k}x\to v$, where $u,v\in X\setminus\Lambda$. Since $X\setminus\Lambda$ is open, one can find $k_0$  such that $T_1^{n_k}x, T_2^{n_k}x\in X\setminus\Lambda$ for all $k\geq k_0$. 
Then the set $K:= \{ T_1^{n_k}x;\; k\geq k_0\}\cup \{ T_2^{n_k}x;\; k\geq k_0\} \cup\{ u, v\}$ is a compact subset of $X\setminus\Lambda$ and $x\in T_1^{-n_k}(K)\cap T_2^{-n_k}(K)$ for all $k\geq k_0$, so that (1) is not satisfied.
\epf

\
Here is a first application of Theorem \ref{Th1}.
\bco\label{ab} Let $X$ be a Polish topological vector space, and let $T\in\mathcal L(X)$. Let also $a,b\in\K$, and assume that $aT$ and $bT$ admit non-trivial invariant measures. Then, $aT$ and $bT$ are orthogonal if and only if $\vert a\vert\neq \vert b\vert$. More precisely: $aT$ and $bT$ are orthogonal if $\vert a\vert\neq \vert b\vert$; and they share a non-trivial invariant measure if  $\vert a\vert=\vert b\vert$.
\eco
\bpf To prove that $aT$ and $bT$ are orthogonal if $\vert a\vert\neq\vert b\vert$, we use Theorem \ref{Th1} and Remark \ref{unpeuutile} with $\Lambda=\{ 0\}$. Assume that there exists $x\in X$ and an increasing sequence of integers $(n_k)$ such that $a^{n_k}T^{n_k}x\to u\in X$ and $b^{n_k}T^{n_k}x\to v\in X$. We must show that $u=0$ or $v=0$; and this follows from the continuity of the map 
$(\lambda, z)\mapsto \lambda z$: indeed, assuming for example that $\vert a\vert >\vert b\vert$ and writing 
$b^{n_k}T^{n_k}x= (b/a)^{n_k} \, a^{n_k}T^{n_k} x$, we see that $b^{n_k}T^{n_k}x\to 0$.

\smallskip The converse will follow from Proposition \ref{symm} below.
\epf

\begin{remark} If $B_\w$ is a weighted shift acting on $X=\ell_p(\Z_+)$ and $a,b\in \K$ are such that $\vert a\vert =\vert b\vert=:r$ with $\sum_{n=1}^\infty \frac{1}{(r^n \vert w_1\cdots w_n\vert)^p}<\infty$, then the operators $aB_{\pmb w}$ and $bB_{\pmb w}$ admit equivalent ergodic Gaussian measures with full support. This will follow from Theorem \ref{gauss} below.
\end{remark}

\smallskip Corollary \ref{ab} leads to naive speculations regarding orthogonality of weighted shifts. Recall that two operators $T_1,T_2$ acting on $X$ are said to be \emph{similar} if there exists an invertible operator $J$ such that $T_1=JT_2J^{-1}$. By \cite{Sh}*[p. 54, Th. 2(a)], two weighted shifts $B_{\pmb u}$ and $B_{\pmb v}$ acting on $\ell_p$ are similar if and only if 

\begin{equation}\label{sim}
\varliminf\; \left\vert\frac{u_1\cdots u_n}{v_1\cdots v_n}\right\vert>0\qquad{\rm and}\qquad \varlimsup\; \left\vert\frac{u_1\cdots u_n}{v_1\cdots v_n}\right\vert<\infty.
\end{equation}

When $\pmb u$ and $\pmb v$ have the form $\pmb u=a \pmb w$ and $\pmb v=b\pmb w$ for some weight sequence $\pmb w$ and $a,b\in\K$, then $\frac{u_1\cdots u_n}{v_1\cdots v_n}=\left(\frac{a}b\right)^n$, and hence (\ref{sim}) holds true if and only if $\vert a\vert=\vert b\vert$. So it is tempting to ``conjecture'' that two weighted shifts $B_{\pmb u}$ and $B_{\pmb v}$ are orthogonal if and only if they are not similar. This is however not true, as will be seen in Section \ref{Sec4} below. Nevertheless, we now present two results showing that this is not \emph{that far} from being true. 

\smallskip
In what follows we denote by $e_j^*$, $j\geq 0$ the coordinate functionals on $X=\ell_p(\Z_+)$. And if $x\in X$, we write $x=\sum_{j=0}^\infty x_j e_j$ where $x_j=\langle e_j^*,x\rangle$.

\bpr\label{presque} Let $B_\u$ and $B_\v$ be two weighted shifts acting on $X=\ell_p$. Assume that $B_\u$ and $B_\v$ are not similar, \textit{i.e.} that
\begin{equation}\label{liminf}
\hbox{either}\quad \underline\lim\, \left\vert \frac{u_1\cdots u_n}{v_1\cdots v_n}\right\vert=0
\qquad\hbox{or}\quad \overline\lim\, \left\vert \frac{u_1\cdots u_n}{v_1\cdots v_n}\right\vert=\infty.
\end{equation}

Then, $B_\u$ and $B_\v$ are orthogonal with respect to measures not charging $\Lambda_{e_0^*}:=\ker(e_0^*)$.
\epr
\bpf We apply Theorem \ref{Th1} with $\Lambda:=\Lambda_{e_0^*}$. So let $K$ be a compact subset of $X\setminus \ker(e_0^*)$. Choose some constants $M<\infty$ and $\gamma >0$ such that $\Vert z\Vert\leq M$ and $\vert \langle e_0^*,z\rangle\vert\geq \gamma$ for all $z\in K$. 
If $x=\sum_{j=0}^\infty x_j e_j\in X$ then, for every $n\geq 1$, we have 
\[ \langle e_0^*, B_\u^nx\rangle =\uprod\, x_n\qquad{\rm and}\qquad \langle e_0^*, B_\v^nx\rangle =\vprod\, x_n;\]
so we see that
\[ B_\u^{-n}(K)\subseteq \bigl\{ x\in X;\; \vert x_n\vert\cdot\vert u_1\cdots u_n\vert\geq \gamma\bigr\}\]
and
\[ B_\v^{-n}(K)\subseteq \bigl\{ x\in X;\; \vert x_n\vert\cdot\vert v_1\cdots v_n\vert\leq M\bigr\}.\]

Hence, the following implication holds :
\[ B_\u^{-n}(K)\cap B_\v^{-n}(K)\neq\emptyset\implies \left\vert \frac{u_1\cdots u_n}{v_1\cdots v_n}\right\vert\geq \frac\gamma{M} \cdot\]

It follows that if $\underline\lim\, \left\vert \frac{u_1\cdots u_n}{v_1\cdots v_n}\right\vert=0$, then $B_\u^{-n}(K)\cap B_\u^{-n}(K)=\emptyset$ for infinitely many $n$; and symmetrically, if $\overline\lim\, \left\vert \frac{u_1\cdots u_n}{v_1\cdots v_n}\right\vert=\infty$ then $B_\v^{-n}(K)\cap B_\u^{-n}(K)=\emptyset$ for infinitely many $n$. By Theorem \ref{Th1}, this concludes the proof.
\epf

\begin{remark} One may observe that for any weight sequence $\pmb w$ and every $n\geq 1$, we have $B_\w^{-n}(\Lambda_{e_0^*})=\Lambda_{e_n^*}$. Hence, a measure $m\in \mathcal P_\w(X)$ does not charge $\Lambda_{e_0^*}$ if and only if it does not charge $\Lambda_{e_n^*}$. Therefore, $m$ does not charge $\Lambda_{e_0^*}$ if and only if it is supported on the set $\{ x\in X;\; \langle e_n^*,x\rangle\neq0\;\,\hbox{for all $n\geq 0$}\}$.
\end{remark}

\smallskip The next theorem shows that ``full'' orthogonality of $B_\u$ and $B_\v$ can be deduced from stronger assumptions on the weights $\u$ and $\v$.
\bth\label{presquebis} Let $B_\u$ and $B_\v$ be two weighted shifts acting on $X=\ell_p$. Assume that 
either
\begin{equation}\label{liminfbis} 
\varliminf_{n\to\infty} \;\max_{0\leq d\leq N} \,\left\vert\frac{u_1\cdots u_{n+d}}{v_1\cdots v_{n+d}}\right\vert=0 \qquad\hbox{for all $N\geq 0$}
\end{equation}
or
\begin{equation}\label{limsupbis} 
\varlimsup_{n\to\infty} \;\min_{0\leq d\leq N} \,\left\vert\frac{u_1\cdots u_{n+d}}{v_1\cdots v_{n+d}}\right\vert=\infty \qquad\hbox{for all $N\geq 0$}.
\end{equation}

\smallskip
Then, $B_\u$ and $B_\v$ are orthogonal.
\eth
\bpf We apply Theorem \ref{Th1} once again, with $\Lambda:=\{ 0\}$. So, let $K$ be a compact subset of $X\setminus\{ 0\}$.

\smallskip We first note that there exists an integer $N$ such that $\max\limits_{0\leq d\leq N} \vert z_d\vert>0$ for every $z\in K$. Indeed, the compact set $K$ is covered by the open sets $O_d:=\{ z\in X;\; \vert z_d\vert>0\}$. 

\smallskip Choose $\gamma>0$ such that $\max\limits_{0\leq d\leq N} \vert z_d\vert\geq \gamma$ for every $z\in K$, and $M<\infty$ such that $\Vert z\Vert\leq M$ for every $z\in K$. As in the proof of Proposition \ref{presque}, we get that for every $n\geq 1$,
\[ B_\u^{-n}(K)\subseteq\left\{ x\in X;\; \max_{0\leq d\leq N}  \vert u_{1+d}\cdots u_{n+d}\vert \cdot\vert x_{n+d}\vert\geq\gamma\right\}\]
and
\[ B_\v^{-n}(K)\subseteq\left\{ x\in X;\; \max_{0\leq d\leq N} \vert v_{1+d}\cdots v_{n+d}\vert \cdot\vert x_{n+d}\vert\leq M\right\};\]
so that
\[ B_\u^{-n}(K)\cap B_\v^{-n}(K)\neq\emptyset\implies  \max_{0\leq d\leq N} \left\vert\frac{u_{1+d}\cdots u_{n+d}}{v_{1+d}\cdots v_{n+d}}\right\vert\geq \frac\gamma{M}\cdot \]

Now, we have $\max_{0\leq d\leq N} \left\vert\frac{u_{1+d}\cdots u_{n+d}}{v_{1+d}\cdots v_{n+d}}\right\vert\leq C_N\, \max_{0\leq d\leq N} \left\vert\frac{u_{1}\cdots u_{n+d}}{v_{1}\cdots v_{n+d}}\right\vert$ for some constant $C_N$. Hence, if (\ref{liminfbis}) is satisfied, then one can find an infinite set $I\subseteq\N$ such that \[ \max\limits_{0\leq d\leq N} \left\vert\frac{u_{1+d}\cdots u_{n+d}}{v_{1+d}\cdots v_{n+d}}\right\vert\to 0\qquad \hbox{as}\quad n\to\infty\;\;\hbox{along $I$};\] and it follows that $B_\u^{-n}(K)\cap B_\v^{-n}(K)=\emptyset$ for all 
$n\in I$ sufficiently large. Symmetrically, if (\ref{limsupbis}) is satisfied, then one can find an infinite set $I\subseteq\N$ such that  $B_\v^{-n}(K)\cap B_\u^{-n}(K)=\emptyset$ for all 
$n\in I$ sufficiently large. So Theorem \ref{Th1} applies.
\epf

\bco\label{bdedbelow} Assume that the weight sequences $\u$ ad $\v$ are bounded below, \textit{i.e.} $\inf_{n\geq 1} \vert u_n\vert>0$ and $\inf_{n\geq 1} \vert u_n\vert>0$. If $B_\u$ and $B_\v$ are not similar, then they are orthogonal. 
\eco
\bpf The assumption on $\u$ and $\v$ implies that for every $N\geq 0$, there are constants $c_N>0$ and $C_N<\infty$ such that 
$\max_{0\leq d\leq N} \,\left\vert\frac{u_1\cdots u_{n+d}}{v_1\cdots v_{n+d}}\right\vert\le C_N\, \left\vert\frac{\uprod}{\vprod}\right\vert$ and $\min_{0\leq d\leq N} \,\left\vert\frac{u_1\cdots u_{n+d}}{v_1\cdots v_{n+d}}\right\vert\geq c_N\, \left\vert\frac{\uprod}{\vprod}\right\vert$. Hence, if $B_\u$ and $B_\v$ are not similar then (\ref{liminfbis}) or (\ref{limsupbis}) is satisfied.
\epf

\smallskip
\begin{remark} Theorem \ref{presquebis} can be reformulated in the following way: if 
\be
\item[-] either, for every $\varepsilon >0$, the set $A_{\u,\v,\varepsilon}:=\left\{ n\in\Z_+;\; \left\vert\frac{u_1\cdots u_n}{v_1\cdots v_n}\right\vert<\varepsilon\right\}$ contains arbitrarily long intervals,
\item[-] or, for every $M<\infty$, the set $A'_{\u,\v,M}:=\left\{ n\in\Z_+;\; \left\vert\frac{u_1\cdots u_n}{v_1\cdots v_n}\right\vert> M\right\}$ contains arbitrarily long intervals,
\ee

then $B_\u$ and $B_\v$ are orthogonal.

\end{remark}

\section{Non-orthogonality \textit{via} periodic points}\label{Sec4} In this section, our aim is to  present some examples of non-orthogonal weighted shifts using measures supported by periodic orbits. We first note the following easy fact.
\begin{fact}\label{per} Let $X$ be a Polish topological vector space, and let $T_1, T_2\in\mathcal L(X)$. If $T_1$ and $T_2$ share a non-zero periodic point, then they are not orthogonal.
\end{fact}
\bpf Let $x\neq 0$ be a common periodic point for $T_1$ and $T_2$; so $T^{d_1}x=x=T^{d_2}x$ for some $d_1,d_2\in\N$. For $i=1,2$, the measure $m_i:=\frac1{d_i}\sum_{n=0}^{d_i-1} \delta_{T_i^nx}$ is $T_i\,$-$\,$invariant, and $m_i(\{ 0\}) =0$ since $x\neq 0$. Moreover, $m_1$ and $m_2$ are not orthogonal since 
$m_1(\{ x\})>0$ and $m_2(\{ x\})>0$.
\epf

\bth\label{pershift} Let $B_\u$ and $B_\v$ be two weighted shifts on $X=\ell_p(\Z_+)$. Assume that $\sum_{n=1}^\infty \frac1{\vert u_1\cdots u_n\vert^p}<\infty$ and $\sum_{n=1}^\infty \frac1{\vert v_1\cdots v_n\vert^p}<\infty$. The following are equivalent.
\be
\item[\rm (1)] $B_\u$ and $B_\v$ share a non-zero periodic point;
\item[\rm (2)] there exist $d\in\N$ and $0\leq j\leq d-1$ such that $u_{1+j}\cdots u_{dm+j}=v_{1+j}\cdots v_{dm+j}$ for all $m\geq 1$;
\item[\rm (3)]  there exist $d\in\N$ and $0\leq j\leq d-1$, and a non-zero scalar $C$ such that $u_1\cdots u_{dm+j}=C\, v_1\cdots v_{dm+j}$ for all $m\geq 0$ {\rm (}where an empty product is declared to be equal to $1${\rm )}.
\ee
\eth
\bpf (1)$\implies$(2) Assume that $B_\u$ and $B_\v$ share a non-zero periodic point $x$. Let $d_1$ be the period of $x$ as a periodic point of $T_1$, and let $d_2$ be the period of $x$ as a periodic point of $T_2$. Setting $d:= d_1d_2$, we then have  
$B_\u^d x=x=B_\v^dx$ (and hence 
$B_\u^{md} x=x=B_\v^{md}x$ for all $m\geq 1$.) Writing $x=\sum_{j\geq 0} x_j e_j$, we have 
\[x_{j+md}=\frac1{u_{j+1}\cdots u_{j+md}}\, x_j=\frac1{v_{j+1}\cdots v_{j+md}}\, x_j\qquad\hbox{for all $j\geq 0$}.\]

Since $x\neq 0$, one can find $0\leq j\leq d-1$ such that $x_j\neq 0$, and then $u_{j+md}\cdots u_{j+1}=v_{j+md}\cdots v_{j+1}$ for all $m\geq 1$.

\smallskip
(2)$\implies$(3) This is clear: if (2) is satisfied for some $d$ and $0\leq j\leq d-1$, then (3) is satisfied with the same $d$ and $j$ and $C:=\frac{u_1\cdots u_j}{v_1\cdots v_j}$ (so $C=1$ if $j=0$).

\smallskip
(3)$\implies$(1) Assume that (3) holds true for some $d, j,C$. Since $\sum_{n=1}^\infty \frac1{\vert u_1\cdots u_n\vert^p}<\infty$, the vector 
\[ x:=\sum_{m=0}^\infty \frac{1}{u_1\cdots u_{md+j}}\, e_{md+j}=\frac1C\sum_{m=0}^\infty \frac{1}{v_1\cdots v_{md+j}}\, e_{md+j}\in \ell_p\]
is well-defined, and clearly $B_\u^dx=x=B_\v^dx$. 
Thus, $x$ is a non-zero common periodic point for $B_\u$ and $B_\v$.
\epf

\bco Let $B_\u$ and $B_\v$ be two weighted shifts on $\ell_p(\Z_+)$, with $\sum_{n=1}^\infty \frac1{\vert u_1\cdots u_n\vert^p}<\infty$ and $\sum_{n=1}^\infty \frac1{\vert v_1\cdots v_n\vert^p}<\infty$. Assume that there exist $d\in\N$, $0\leq j\leq d-1$ and a scalar $C\neq 0$ such that $u_1\cdots u_{dm+j}=C\, v_1\cdots v_{dm+j}$ for all $m\geq 0$. Then $B_\u$ and $B_\v$ are not orthogonal. 
\eco

\smallskip
\begin{example} Assume that $u_j=2$ for all $j\geq 1$, $v_1=3$ and $v_j=2$ for all $j\geq 2$. Then $B_\u$ and $B_\v$ acting on $\ell_p$ share a non-zero fixed point, and hence they are not orthogonal.
\end{example}

\smallskip In this example, the operators $B_\u$ and $B_\v$ are obviously similar. We now give an example of non-orthogonal weighted shifts which are not similar. 
\begin{example}\label{non} Let $(r_k)_{k\geq 1}$ be an increasing sequence of integers, and let $\u$ and $\v$ be the two weight sequences defined as follows:
\[ u_j:=\left\{ \begin{array}{lll}
2&\hbox{if}&\hbox{$j\not\in\{ 5r_k+1, 5r_k+4;\; k\geq 1\}$}\\
\frac1k&\hbox{if}& \hbox{$j=5r_k+1$ or $j=5r_k+4\;$ for some $k\geq 1$}
\end{array}\right.\]
and
\[ v_j:=\left\{ \begin{array}{lll}
2&\hbox{if}&\hbox{$j\not\in\{ 5r_k+2, 5r_k+3;\; k\geq 1\}$}\\
\frac1k&\hbox{if}& \hbox{$j=5r_k+2$ or $j=5r_k+3\;$ for some $k\geq 1$.}
\end{array}\right.\]

\smallskip
If $(r_k)$ is sufficiently fast increasing, then $B_\u$ and $B_\v$ acting on $\ell_p$ share a non-zero periodic point and hence they are not orthogonal. However, we have
\[ \varliminf \;\frac{u_1\cdots u_n}{v_1\cdots v_n}=0\qquad{\rm and}\qquad \varlimsup\; \frac{u_1\cdots u_n}{v_1\cdots v_n}=\infty,\]
so that, in particular, $B_\u$ and $B_\v$ are not similar.
\end{example}
\bpf It is clear that if $(r_k)$ is sufficiently fast increasing, then $\sum_{n=0}^\infty \frac1{\vert u_1\cdots u_n\vert^p}<\infty$ and $\sum_{n=0}^\infty \frac1{\vert v_1\cdots v_n\vert^p}<\infty$.

Observe that if $m$ is a positive integer and if we consider the largest $k$ such that $r_k\leq m$, then either $r_k=m$ or $5r_k+4<5m$. Since  $u_1\cdots u_{5r_k}=v_1\cdots v_{5r_k}$ by definition of the weights, it follows  that 
\[ u_1\cdots u_{5m}=v_1\cdots v_{5m}\qquad\hbox{for all $m\geq 1$}.\]
Hence, $B_\u$ and $B_\v$ share a non-zero periodic point by Theorem \ref{pershift}.  

However, for any $k\geq 1$ we have 
\[ \frac{u_1\cdots u_{5r_k+1}}{v_1\cdots v_{5r_k+1}}=\frac1{2k}\qquad{\rm and}\qquad \frac{u_1\cdots u_{5r_k+3}}{v_1\cdots v_{5r_k+3}}=2k,\]
so that $\varliminf \;\frac{u_1\cdots u_n}{v_1\cdots v_n}=0$ and $\varlimsup\; \frac{u_1\cdots u_n}{v_1\cdots v_n}=\infty$.

\epf

\begin{remark} In the above example, the weights are not bounded below. This is not accidental: indeed, if two weighted shifts $B_\u$ and $B_\v$ share a non-zero periodic point and if 
$\inf_{n\geq 1} \vert u_n\vert >0$ and $\inf_{n\geq 1} \vert v_n\vert >0$, then $B_\u$ and $B_\v$ are necessarily similar by Corollary \ref{bdedbelow}. (Alternatively, this follows easily from Theorem \ref{pershift}.)  
\end{remark}


\section{Non-orthogonality \textit{via} product measures}\label{Sec5} 

\subsection{Invariant product measures for weighted shifts}
The question we consider in this section is the following: given two weighted shifts $B_\u$ and $B_\v$ on $X=\ell_p(\Z_+)$, when is it possible to find invariant measures $m_\u$ and $m_\v$ (for $B_\u$ and $B_\v$ respectively) not charging $\{ 0\}$ which are not orthogonal and are also \emph{product measures}?

\smallskip
The terminology requires some explanation since product measures should be defined on products of measurable spaces and $X=\ell_p$ is not  such a product space. However, $\ell_p$ is contained in the product space $\Omega:=\K^{\Z_+}$. It is well-known that the Borel $\sigma$-algebra of $\Omega$ (induced by the product topology) is identical with the product $\sigma$-algebra $\otimes_{n\geq 0} \mathcal B_\K$. Moreover, $\ell_p$ is a Borel subset of $\Omega$ when 
$\Omega$ is endowed with the product topology, and a subset of $\ell_p$ is Borel in $(\ell_p,\Vert\,\cdot\,\Vert_p)$ if and only if is Borel in $\Omega$.  It follows that if $\mu$ is a Borel probability measure on $\Omega$, then the restriction of 
$\mu$ to $\ell_p$ is a Borel measure on $\ell_p$ endowed with its usual topology (but not a probability measure unless $\mu(\ell_p)=1$); and conversely, any Borel probability measure on $\ell_p$ can be considered as a Borel probability measure on $\Omega$ (supported on $\ell_p$). We will say that a Borel probability measure $m$ on $\ell_p$ is a \emph{product measure on $\ell_p$} if $m$ is the restriction to $\ell_p$ of a (probability) product measure $\mu$ on $\Omega$, \textit{i.e.} $\mu=\bigotimes_{n\geq 0} \mu_n$ where each $\mu_n$ is a Borel probability measure on $\K$,  such that $\mu(\ell_p)=1$. In this case, the measure $\mu$ is uniquely determined by $m$ since 
$\mu(B)=m(B\cap \ell_p)$ for every Borel set $B\subseteq \Omega$, so we identify $m$ and $\mu$ and simply write $m=\otimes_{n\geq 0} \mu_n$.

\smallskip
The following  essentially obvious remark will be used repeatedly.
\begin{fact} Let $X$ be a Borel subset of $\Omega= \K^{\Z_+}$, and let $\mu$ be a Borel (probability) measure on $\Omega$ such that $\mu(X)=1$. Set $m:=\mu_{| X}$. Let also $T:\Omega\to\Omega$ be a Borel map such that $T(X)\subseteq X$.
Then $\mu$ is $T\,$-$\,$invariant if and only if $m$ is 
$(T_{| X})\,$-$\,$invariant. 
\end{fact}
\bpf The measure $\mu$ is $T\,$-$\,$invariant if and only if $\mu(T^{-1}(B))=\mu(B)$ for every Borel set $B\subseteq\Omega$. Now, we have $\mu(B)=m(X\cap B)$ and $\mu(T^{-1}(B))=m(X\cap T^{-1}(B))=m\bigl( (T_{| X})^{-1}(X\cap B)\bigr)$ because $T(X)\subseteq X$. So $\mu$ is $T\,$-$\,$invariant 
if and only if $m(X\cap B)=m\bigl( (T_{| X})^{-1}(X\cap B)\bigr)$ for every Borel set $B\subseteq\Omega$, which means exactly that $m$ is $(T_{| X})\,$-$\,$invariant.
\epf

\smallskip From this observation, it follows that if $m$ is a product measure on $X=\ell_p$ and if we denote by $\mu$ the measure $m$ considered as a measure on $\Omega=\K^{\Z_+}$, then $m$ is invariant under some weighted shift $B_\w$ if and only if $\mu$ is invariant under the natural extension of $B_\w$ to $\Omega$ (defined by the same formula and also denoted by $B_\w$). The next lemma says precisely when this happens.
\blm\label{invar} Let $\mu=\otimes_{n\geq 0} \mu_n$ be a product measure on $\Omega=\K^{\Z_+}$, and let $\w$ be a weight sequence. Then, $\mu$ is $B_\w\,$-$\,$invariant if and only if, for each $n\geq 1$, the measure $\mu_n$ is the image of $\mu_0$ under the map $t\mapsto \frac1{w_1\cdots w_n}\, t$, \textit{i.e.} $\mu_n(A)=\mu_0(w_1\cdots w_nA)$ for every Borel set $A\subseteq \K$.
\elm
\bpf 
For any Borel sets $A_0\,\dots ,A_N\subseteq\K$,  denote by $[A_0,\dots ,A_N]\subseteq \Omega$ the ``cylinder set'' defined as follows:
\[ [A_0,\dots ,A_N]=\bigl\{ t=(t_j)_{j\geq 0}\in\Omega;\; t_j\in A_j\quad\hbox{for $j=0,\dots ,N$}\bigr\}.\]

Since the cylinder sets generate the Borel $\sigma$-algebra of $\Omega$, the measure $\mu$ is $B_\w\,$-$\,$invariant if and only if 
\[  \mu\bigl( B_\w^{-1}([A_0,\dots ,A_N])\bigr)=\mu\bigl([A_0,\dots ,A_N]\bigr)\]
for all Borel sets $A_0,\dots ,A_N\subseteq \K$. Now, by definition of $B_\w$ we have 
\[ B_\w^{-1}([A_0,\dots ,A_N]\bigr)=\bigl\{ t=(t_j)_{j\geq 0}\in\Omega;\; w_{j+1}t_{j+1}\in A_j\quad\hbox{for $j=0,\dots ,N$}\bigr\},\]
so that 
\[ B_\w^{-1}([A_0,\dots ,A_N]\bigr)=\left[\K, (1/{w_1}) A_0,\dots ,(1/{w_{N+1}}) A_N\right].\]
Since $\mu=\otimes_{j\geq 0} \mu_j$, it follows that $\mu$ is $B_\w\,$-$\,$invariant if and only if 
\[ \prod_{j=0}^N \mu_j(A_j)= \prod_{j=0}^{N} \mu_{j+1}\bigl((1/w_{j+1}) A_{j}\bigr)\]
for every $N\geq 0$ and all Borel sets $A_0,\dots ,A_N\subseteq \K$. This is clearly equivalent to the fact that
\[ \mu_{j+1}((1/w_{j+1})A)=\mu_{j}(A)\qquad\hbox{for all $j\geq 0$ and every Borel set $A\subseteq\K$},\]
which proves the lemma.
\epf

\bco  Let $\mu=\otimes_{n\geq 0} \mu_n$ be a product measure on $\Omega=\K^{\Z_+}$, and let $\w$ be a weight sequence. Assume that $\mu_0$ has a density $\mathtt p$ with respect to Lebesgue measure on $\K$. 
 Then $\mu$ is $B_\w\,$-$\,$invariant if and only if, for each $n\geq 1$, the measure $\mu_n$ has a density $\mathtt p_n$ given by $\mathtt p_n(t):=  \vert w_1\cdots w_n\vert^d\,  \mathtt p(w_1\cdots w_n t)$, where $d=1$ if $\K=\R$ and $d=2$ if $\K=\C$.
\eco
\bpf This follows immediately from Lemma \ref{invar}
\epf

\bco\label{charge} Let $\mu=\otimes_{n\geq 0} \mu_n$ be a product measure on $\Omega=\K^{\Z_+}$, and let $\w$ be a weight sequence such that $\sum_{n=0}^\infty \frac1{\vert w_1\cdots w_n\vert^p}<\infty$. If $\mu$ is $B_\w\,$-$\,$invariant and if the measure $\mu_0$ is such that $\int_\K \vert t\vert^p d\mu_0(t)<\infty$, then $\mu(\ell_p)=1$.
\eco
\bpf By Lemma \ref{invar}, we have 
\begin{align*}
\int_\Omega \left(\sum_{n=0}^\infty \vert t_n\vert^p\right) d\mu(t)&=\sum_{n=0}^\infty \int_\K \vert t_n\vert^p d\mu_n(t_n)\\
&
=\left(\int_\K \vert s\vert^p d\mu_0(s)\right)\times\left(1+\sum_{n=1}^\infty \frac1{\vert w_1\cdots w_n\vert^p}\right) <\infty.
\end{align*}

In particular, we see that $\sum_{n=0}^\infty \vert t_n\vert^p<\infty$ for $\mu\,$-$\,$almost every $t\in\Omega$, \textit{i.e.} $\mu(\ell_p)=1$.
\epf

\smallskip To conclude this section, we point out the following general fact.
\begin{fact}\label{ergprod} Let $B_\w$ be a weighted shift on $\ell_p$ where $\w$ satisfies $\sum_{n=1}^\infty \frac1{\vert\wprod\vert^p}<\infty$, and let $m=\otimes_{n\geq 0} \mu_n$ be a $B_\w\,$-$\,$invariant product measure on $\ell_p$. If the measure $\mu_0$ has full support, then $m$ is an ergodic measure with full support for $B_\w$.
\end{fact}
\bpf Let $(\xi_n)_{n\geq 0}$ be a sequence of independent $\K\,$-$\,$valued random variables with law $\mu_0$. As a measure on $\K^{\Z_+}$, $m$ is the distribution of the $\K^{\Z_+}$-$\,$valued random variable 
$\xi:=\xi_0e_0+\sum_{n=1}^\infty \frac{\xi_n}{\wprod}\, e_n$. Since $m(\ell_p)=1$, we have $\xi\in\ell_p$ almost surely, which means that the series defining $\xi$ is almost surely convergent with respect to the $\ell_p\,$-$\,$norm. Moreover, since $\mu_0$ has full support, the measure $m$ has full support. By \cite{Kevin}*[Proposition 2.5], it follows that $m$ is an ergodic measure (in fact, a strongly mixing measure) for $B_\w$.
\epf

\subsection{Non-orthogonality of product measures} A basic tool in the study of product measures is a famous classical result of Kakutani \cite{Kak} giving a necessary and sufficient condition for the equivalence of two infinite product probability measures $\nu=\otimes_{n\geq 0}\nu_n$ and $\nu'=\otimes_{n\geq 0}\nu'_n$ defined on a product of measurable spaces $(\Omega,\mathcal B)=\otimes_{n\geq 0} (\Omega_n,\mathcal B_n)$. 

\smallskip
Let us first recall the definition of the so-called \emph{Hellinger integral} $H(\alpha, \beta)$ of two probability measures $\alpha,\beta$ on some measurable space $(T,\mathcal T)$:
\[ H(\alpha,\beta)=\int_T \sqrt{\frac{d\alpha}{d\tau}}\,\sqrt{\frac{d\beta}{d\tau}}\, d\tau,\]
where $\tau$ is any sigma-finite positive measure on $(T,\mathcal T)$ such that $\alpha$ and $\beta$ are absolutely continuous with respect to $\tau$. If we write 
$\alpha =\frac{d\alpha}{d\beta} \, \beta + \nu$ the Lebesgue-Nikodym decomposition of $\alpha$ with respect to $\beta$, then $\left(\frac{d\alpha}{d\tau}-\frac{d\alpha}{d\beta}\, \frac{d\beta}{d\tau}\right) \frac{d\beta}{d\tau}=0$ $\tau\,$-$\,$a.e. since $\nu\perp \beta$, so $\frac{d\alpha}{d\tau}\, \frac{d\beta}{d\tau}= \frac{d\alpha}{d\beta} \left(\frac{d\beta}{d\tau}\right)^2$ and hence 
\[ H(\alpha,\beta)=\int_T \sqrt{\frac{d\alpha}{d\beta} }\, d\beta;\]
which shows that $H(\alpha,\beta)$ is indeed independent of the choice of $\tau$. Note that we always have $0\leq H(\alpha,\beta)\leq 1$ by the Cauchy-Schwarz inequality; and that $H(\alpha,\beta)=1$ if and only if $\alpha=\beta$, whereas $H(\alpha,\beta)=0$ if and only if $\alpha\perp\beta$. 

\smallskip
Kakutani's Theorem as we will need it can now be stated as follows.  
\bth Let $\nu=\otimes_{n\geq 0}\nu_n$ and $\nu'=\otimes_{n\geq 0}\nu'_n$ be two product probability measures on $(\Omega,\mathcal B)=\otimes_{n\geq 0} (\Omega_n,\mathcal B_n)$. The measures $\nu$ and $\nu'$ are non-orthogonal if and only if 
\[ \prod_{n=0}^\infty H(\nu_n,\nu'_n)>0.\]
Moreover, under the assumption that $\nu_n\sim \nu'_n$ for each $n\geq 0$, the measures $\nu$ and $\nu'$ are either orthogonal or equivalent.
\eth 

\smallskip
Note that this is not exactly  what is proved in \cite{Kak}: the main result of \cite{Kak} states that under the assumption that $\nu_n\sim \nu'_n$ for all $n$, the measures $\nu$ and $\nu'$ are equivalent if $\prod_{n=0}^\infty H(\nu_n,\nu'_n)>0$ and orthogonal otherwise. However, the above version of the theorem is certainly well-known (it is stated for example in \cite{Shepp}), and it can be obtained by slight modifications of Kakutani's original proof; see Section \ref{add} for more details.

\smallskip
\subsection{Gaussian product measures}
A particularly interesting class of product measures is that of \emph{Gaussian product measures}. Our definition will be very restrictive:  
we will say that a product measure $\mu=\otimes_{n\geq 0} \mu_n$ on $\Omega=\K^{\Z_+}$ is a Gaussian product measure if for each $n\geq 0$, the measure $\mu_n$ is a  Gaussian measure on $\K$ of the form $\mathcal N(0,\sigma_n^2)$ for some $\sigma_n>0$, \textit{i.e.} $\mu_n$ is a centered Gaussian measure with  covariance matrix $\sigma_n^2I_\K$. Accordingly, we say that a measure $m$ on $X=\ell_p$ is a  Gaussian product measure if $m$ is the restriction to $\ell_p$ of a  Gaussian product measure $\mu$ on $\Omega$ such that $\mu(\ell_p)=1$.

\smallskip By Fernique's integrability theorem (see e.g. \cite{Bog}), if $m$ is a Gaussian measure on $X=\ell_p$, then $m$ admits moments of all orders and in particular $\int_X \Vert x\Vert^p dm<\infty$. It follows that if $\mu=\otimes_{n\geq 0} \mu_n$ is a  Gaussian product measure on $\Omega=\K^{\Z_+}$, then $\mu$ is supported on $\ell_p$ if and only if $\int_\Omega \bigl(\sum_{n=0}^\infty \vert t_n\vert^p\bigr) \, d\mu(t)<\infty$. Moreover, if $\mu$ is invariant under some weighted shift $B_\w$, then $\int_\Omega \vert t_n\vert^p d\mu(t)=\frac1{\vert w_1\cdots w_n\vert^p}\, \int_\K \vert t\vert^p d\mu_0(t)$ for all $n\geq 1$ by Lemma \ref{invar}. Altogether, we see that a weighted shift $B_\w$ acting on $\ell_p$ admits an invariant  Gaussian product measure if and only if $\sum_{n=0}^\infty\frac1{\vert w_1\cdots w_n\vert^p}<\infty$, and that if this holds, then any $B_\w\,$-$\,
$invariant Gaussian product measure on $\K^{\Z_+}$ is supported on $\ell_p$. (For the ``only if'' part, one could also have used Proposition \ref{existe} rather than Fernique's integrability theorem.) We refer the reader to \cite{BG1}, \cite{BG2}, \cite{BM3} and \cite{BM}*{Chapter 5} for an in-depth study of invariant Gaussian measures for operators on Hilbert or Banach spaces. We just point out here that if $\w$ is a weight sequence such that $\sum_{n=0}^\infty\frac1{\vert w_1\cdots w_n\vert^p}<\infty$ and if $m$ is a
$B_\w\,$-$\,$invariant Gaussian product measure on $\ell_p$, then $m$ is an ergodic measure with full support for $B_\w$ (see Fact \ref{ergprod}).

\smallskip
The next theorem will provide a necessary and sufficient condition for two weighted backward shifts $B_\u$ and $B_\v$ on $X=\ell_p$ to admit equivalent  invariant Gaussian product measures.
\bth\label{gauss} Let $\u$ and $\v$ be two weight sequences, and let $\mu_\u=\otimes_{n\geq 0} \mu_{\u,n}$ and $\mu_\v=\otimes_{n\geq 0} \mu_{\v,n}$ be two Gaussian product measures on $\K^{\Z_+}$ invariant under $B_\u$ and $B_\v$ respectively, with $\mu_{\u,0}=\mathcal N(0,\sigma^2)$ and $\mu_{\v,0}=\mathcal N(0,\sigma'^2)$.  Then $\mu_{\u}$ and $\mu_{\v}$ are either equivalent or orthogonal, and they are equivalent if and only if 
\begin{equation}\label{gaussequ} \sum_{n=1}^\infty \left( 1-{\frac{\sigma'}{\sigma}} \,{\frac{\vert\uprod\vert}{\vert\vprod\vert}}\, \right)^2<\infty.
\end{equation}
\eth

From this result, we immediately deduce
\bco\label{gauss'} Let $B_\u$ and $B_\v$ be two weighted shifts on $X=\ell_p$, where the weight sequences $\u$ and $\v$ are such that   $\sum_{n=0}^\infty\frac1{\vert u_1\cdots u_n\vert^p}<\infty$ and  $\sum_{n=0}^\infty\frac1{\vert v_1\cdots v_n\vert^p}<\infty$. The following are equivalent.
\be
\item[\rm (a)] $B_\u$ and $B_\v$ admit equivalent  invariant Gaussian product measures.
\item[\rm (b)] There exists a constant $\kappa>0$ such that $\sum\limits_{n=1}^\infty \left( 1-\kappa\,{ \frac{\vert \uprod\vert }{\vert\vprod\vert}}\,\right)^2<\infty$.
\ee
\eco

\smallskip Here is another consequence of Theorem \ref{gauss}. 
\bco\label{debile} Let $\pmb\Lambda$ be a countable family of weight sequences. Assume that $\sum\limits_{n=1}^\infty \frac1{\vert\wprod\vert^p}<\infty$ for all $\w\in\pmb\Lambda$ and $\sum\limits_{n=1}^\infty \left(1-{ \frac{\vert \uprod\vert }{\vert\vprod\vert}}\,\right)^2<\infty$ for any $\mathbf u,\mathbf v\in \pmb\Lambda$. Then, there exists a Gaussian measure with full support $m$ on $X=\ell_p$ with the following property: given any sequence of Borel sets $(A_i)_{i\geq 0}\subseteq X$ with $m(A_i)>0$, there exists $x\in X$ such that for every $\w\in \mathbf \Lambda$ and all $i\geq 0$, the set $\mathcal N_{B_\w}(x, A_i):=\{ n\in \N;\; B_\w^n x\in A_i\}$ has a positive density.
\eco
\bpf For each $\w\in\pmb\Lambda$, let $m_\w=\otimes_{n\geq 0} \mu_{\w,n}$ be the $B_\w\,$-$\,$invariant Gaussian product measure on $\ell_p$ defined by $\mu_{\w,0}:= \mathcal N(0,1)$. This is an ergodic measure with full support for $B_\w$. Since the measures $m_\w$, $\w\in \pmb\Lambda$ are pairwise equivalent by Theorem \ref{gauss} and since $\mathbf\Lambda$ is countable, the result follows from the pointwise ergodic theorem, taking $m:= m_{\w_0}$ for any $\w_0\in\mathbf \Lambda$. 
\epf

\begin{remark}\label{same} Taking as $(A_i)_{i\geq 0}$ a countable basis of open sets for $X$ in Corollary \ref{debile}, we get a vector $x\in X$ which is frequently hypercyclic for all operators $B_\w$, $\w\in\mathbf \Lambda$. However, it turns out that something much stronger holds true: if $\u$ and $\v$ are two weight sequences such that $\frac{\uprod}{\vprod}$ has a non-zero limit as $n\to\infty$, then $B_\u$ and $B_\v$ have in fact \emph{the same} frequently hypercyclic vectors. This result,  which is due to S. Charpentier and the third author, will be proved in Section \ref{add}.
\end{remark}


\smallskip Before proving Theorem \ref{gauss}, we use Corollary \ref{gauss'} to give an example of two non-orthogonal weighted shifts sharing no non-zero periodic points.
\begin{example} Let $\u$ and $\v$ be the weight sequences defined as follows:
\[ \u_n:=2\qquad {\rm and}\qquad \v_n:=2\,\frac{1+\varepsilon_n}{1+\varepsilon_{n-1}},\]
where $\varepsilon_0=0$ and $(\varepsilon_n)_{n\geq 1}$ is a decreasing sequence of positive real numbers such that $\sum_{n=1}^\infty \varepsilon_n^2<\infty$. Then $B_\u$ and $B_\v$ acting on $X=\ell_p$ share no non-zero periodic point, but they admit equivalent  invariant Gaussian product measures (and hence they are not orthogonal).
\end{example}
\bpf We have $\uprod=2^n$ and $\vprod=2^n(1+\varepsilon_n)$ for all $n\geq 1$. In particular, $\sum_{n=1}^\infty \frac1{(\uprod)^p}<\infty$ and  $\sum_{n=1}^\infty \frac1{(\vprod)^p}<\infty$. Since $\frac{\uprod}{\vprod}=\frac1{(1+\varepsilon_n)}$ is increasing, the conditions of Theorem \ref{pershift} cannot be satisfied, so $B_\u$ and $B_\v$ do not share any non-zero periodic point. However 
\[ 1-{\frac{\uprod}{\vprod}}=1-\frac1{1+\varepsilon_n}\sim \varepsilon_n\qquad\hbox{as $n\to\infty$},\]
so $B_\u$ and $B_\v$ admit equivalent  invariant Gaussian product measures by Corollary \ref{gauss'}.
\epf

\smallskip It is now time to prove Theorem  \ref{gauss}. 


\smallskip
\begin{proof}[Proof of Theorem \ref{gauss}] In what follows, we set $d:=1$ if $\K=\R$ and $d:=2$ if $\K=\C$. 

\smallskip By Kakutani's Theorem, we just have to show that $\prod_{n=1}^\infty H(\mu_{\u,n},\mu_{\v,n})>0$ if and only if (\ref{gaussequ}) is satisfied.

\smallskip By Lemma \ref{invar}, the measures $\mu_{\u,n}$ and $\mu_{\v,n}$ are uniquely determined by $\mu_{\u,0}=\mathcal N(0,\sigma^2)$ and $\mu_{\v,0}=\mathcal N(0, \sigma'^2)$. Explicitely,
\[ \mu_{\u,n}=\mathcal N(0,\sigma_n^2)\qquad\hbox{with}\quad \sigma_n^2=\frac{\sigma^2}{\vert\uprod\vert^2}\]
and
\[ \mu_{\v,n}=\mathcal N(0,\sigma_n'^2)\qquad\hbox{with}\quad \sigma_n'^2=\frac{\sigma'^2}{\vert\vprod\vert^2}\cdot\]

 The computation of $H(\mu_{\u,n},\mu_{\v,n})$ now relies on the following fact.

\begin{fact}\label{calcul} If $\mu$ and $\mu'$ are two Gaussian measures on $\K$ of the form $\mu=\mathcal N(0,\sigma^2)$ and $\mu'=\mathcal N(0,\sigma'^2)$, then
\[ H(\mu,\mu')=\left(\frac{2\sigma\sigma'}{\sigma^2+\sigma'^2}\right)^{d/2}=\left(\frac{2(\sigma'/\sigma)}{1+(\sigma'/\sigma)^2}\right)^{d/2}\cdot\]
\end{fact}
\begin{proof}[Proof of Fact \ref{calcul}] The measures $\mu$ and $\mu'$ are absolutely continuous with respect to Lebesgue measure on $\K$ with densities $\frac1{({2\pi}\sigma^2)^{d/2}}e^{-\vert t\vert^2/2\sigma^2}$ and $\frac1{({2\pi}\sigma'^2)^{d/2}}e^{-\vert t\vert^2/2\sigma'^2}$. Hence,
\begin{align*}
H(\mu,\mu')&=\frac1{({2\pi \sigma\sigma')^{d/2}}}\,\int_\R e^{-\frac14\left(\frac1{\sigma^2}+\frac1{\sigma'^2}\right) \vert t\vert^2} dt=\frac1{({2\pi \sigma\sigma')^{d/2}}}\times \left(\frac{4\pi}{\frac1{\sigma^2}+\frac1{\sigma'^2}}\right)^{d/2}\cdot
\end{align*}
\epf

Going back to the measures $m_\u=\otimes_{n\geq 0} \mu_{\u,n}$ and $m_\v=\otimes_{n\geq 0} \mu_{\v,n}$, let us set for each $n\geq 1$:
\[\lambda_n:={\frac{\sigma'_n}{\sigma_n}}={\frac{\sigma'}{\sigma}}\,{\frac{\vert\uprod\vert}{\vert\vprod\vert}}\cdot\]
By Fact \ref{calcul}, we have
\[ H(\mu_{\u,n},\mu_{\v,n})=\left(\frac{2 \, \lambda_n}{{1+\lambda_n^2}}\right)^{d/2}\cdot\]

Since $0\leq H(\mu_{\u,n},\mu_{\v,n})\leq 1$ for all $n\geq 1$, it follows that
\[ \prod_{n=1}^\infty H(\mu_{\u,n},\mu_{\v,n})>0 \quad\hbox{if and only if}\quad  \sum_{n=1}^\infty \left(1-\left(\frac{2 \, \lambda_n}{{1+\lambda_n^2}}\right)^{d/2}\right)<\infty.\]

This can be satisfied only if $\frac{2 \, \lambda_n}{{1+\lambda_n^2}}\to 1$ as $n\to\infty$; and an examination of the function $\Psi(\lambda):=\frac{2 \, \lambda }{{1+\lambda ^2}}$ shows that this is equivalent to the fact that $\lambda_n\to 1$.  In this case, writing $\lambda_n=1+u_n$ and 
using Taylor's formula, we see that 
\[  1-\left(\frac{2 \, \lambda_n}{{1+\lambda_n^2}}\right)^{d/2}\sim \frac{d}{4}\, u_n^2=\frac{d}4\, (1-\lambda_n)^2.\]

So we get that 
\[ \prod_{n=1}^\infty H(\mu_{\u,n},\mu_{\v,n})>0\quad\hbox{if and only if}\quad \sum_{n=1}^\infty (1-\lambda_n)^2<\infty,
\]
which is (\ref{gaussequ}). Theorem \ref{gauss} is proved.
\epf

\smallskip 
\subsection{Products of absolutely continuous measures} It is natural to believe that some condition strictly weaker than (b) of Corollary \ref{gauss'} might still yield the existence of (non-Gaussian) equivalent non-trivial invariant product measures for the backward shifts $B_\u$ and $B_\v$. The next theorem shows that this is in fact not the case, at least for 
product measures whose marginals are absolutely continuous with respect to Lebesgue measure. 

\smallskip We point out that it is possible to deduce parts (1ii) and (2) of this theorem from a Theorem of Shepp \cite{Shepp} concerning (non-)orthogonality of translates of a product measure on $\R^{\Z_+}$, namely \cite{Shepp}*{Theorem 1}. However, we are going to give here a self-contained proof which looks rather different from what is done in \cite{Shepp} (except, of course, for the use of Kakutani's Theorem), and we will indicate in Section \ref{add} how Shepp's result can be used in our context; see Theorem \ref{Shepp}.

\smallskip
For simplicity, we will assume that $\K=\R$ and that all the weights are positive. 

\bth\label{gauss2} Let $\u$ and $\v$ be two positive weight sequences. Assume that the backward shifts $B_\u$ and $B_\v$ acting on the real space $X=\ell_p$ admit invariant product measures $m_\u=\otimes_{n\geq 0} \mu_{\u,n}$ and $m_{\v}=\otimes_{n\geq 0}\mu_{\v,n}$ such that $\mu_{\u,0}$ and $\mu_{\v,0}$ are absolutely continuous with respect to Lebesgue measure on $\R$, $\mu_{\u,0}=\mathtt p(t) dt$ and $\mu_{\v,0}= \mathtt q(t)dt$.
\be
\item[\rm (1)] If $m_\u$ and $m_\v$ are not orthogonal, then: 

\smallskip
\be
\item[\rm (i)] the quotient $\frac{\uprod}{\vprod}$ has a limit $a\in (0,\infty)$ as $n\to \infty$, and $\mathtt q(t)=a \mathtt p(a t)$ almost everywhere;
\item[\rm (ii)] we have $\sum\limits_{n=1}^\infty \left(1-\frac1a\, {\frac{\uprod}{\vprod}}\right)^2<\infty$.
\ee

\smallskip\item[\rm (2)] Assume that the function $f:=\sqrt{\mathtt p}$ is continuous  on $\R\setminus\{ 0\}$ and $\mathcal C^1$-$\,$smooth  except at a finite number of points with  $tf'(t)\in L_2(\R)$, and that $\mathtt q(t)= a \mathtt p(a t)$ for some $a\in (0,\infty)$. Then, $m_\u$ and $m_\v$ are non-orthogonal \emph{if and only if} $\sum\limits_{n=1}^\infty \left(1-\frac1a\, {\frac{\uprod}{\vprod}}\right)^2<\infty$.

\smallskip
\item[\rm (3)] Assume that $f=\sqrt{\mathtt p}$ is $\mathcal C^1$-$\,$smooth  except at a finite number of points with  $tf'(t)\in L_2(\R)$, that $\mathtt p$ has at least one discontinuity point in 
$\R\setminus\{ 0\}$, and that  $\mathtt q(t)= a \mathtt p(a t)$ for some $a\in (0,\infty)$. Then, $m_\u$ and $m_\v$ are non-orthogonal if and only if  
$\sum\limits_{n=1}^\infty \left\vert 1-\frac1a\, {\frac{\uprod}{\vprod}}\right\vert<\infty$.
\ee 
\eth

From part (1) of this theorem and Theorem \ref{gauss}, we immediately deduce

\bco\label{bien} Let $B_\u$ and $B_\v$ be two shifts  with positive weights acting on the real space $\ell_p$. The following are equivalent.
\be
\item[\rm (a)] $B_\u$ and $B_\v$ are not orthogonal with respect to product measures whose marginals are absolutely continuous with respect to Lebesgue measure on $\R$.
\item[\rm (b)] There exists $\kappa>0$ such that $\sum_{n=1}^\infty \left( 1-\kappa\, \frac{\uprod}{\vprod}\right)^2<\infty$.
\item[\rm (c)] $B_\u$ and $B_\v$ admit equivalent invariant Gaussian product measures.
\ee
\eco

We point out that Corollary \ref{bien} is also valid in the complex case and for not necessarily positive weights; see Corollary \ref{bien2}.

\smallskip For the proof of Theorem \ref{gauss2}, we will need the following lemma.

\blm\label{bizarre}
Let $h\in L_2(\R)$ be non-zero and real-valued, and let $\mathcal Ph:\R\to \R$ be the function defined by 
\[ \mathcal Ph(\alpha):=\int_\R h(x) h(x+\alpha)\, dx.\]
\be
\item[\rm (a)] The function $\mathcal Ph$ is continuous and even, with $ \mathcal Ph(\alpha)\leq \Vert h\Vert_2^2$  for avery $\alpha\in\R$ and $\mathcal Ph(0)=\Vert h\Vert_2^2$.
\item[\rm (b)] If $h$ belongs to the Sobolev space $W^{1,2}(\R)$, then $\mathcal Ph$ is $\mathcal C^2$-$\,$smooth on $\R$, with  $(\mathcal Ph)'(0)=0$ and $(\mathcal Ph)''(0)<0$. In particular, $\Vert h\Vert_2^2-\mathcal Ph(\alpha)\sim c\alpha^2$ as $\alpha\to 0$, for some constant $c>0$.
\item[\rm (c)] If $h\not\in W^{1,2}(\R)$, then $(\Vert h\Vert_2^2-\mathcal Ph(\alpha))/\alpha^2\to\infty$ as $\alpha\to 0$.
\item[\rm (d)] Assume that $h$ is $\mathcal C^1$-$\,$smooth except at finitely many points, with $h'\in L_2(\R)$, and that $h$ has at least one discontinuity point. Then $\mathcal Ph$ is left-differentiable and right-differentiable at $0$, with $(\mathcal Ph)'(0^-)>0$ and $(\mathcal Ph)'(0^+)<0$.
\ee
\elm

This lemma has the following immediate consequence.
\bco For any non-zero real-valued $h\in L_2(\R)$, there exists a constant $c>0$ such that $\Vert h\Vert_2^2-\mathcal Ph(\alpha)\geq c\alpha^2$ for $\alpha$ sufficiently close to $0$.
\eco

Note that using Plancherel's formula, it can be seen that this can be stated equivalently as follows: for any $g\in L_2(\R)\setminus\{ 0\}$ with real-valued Fourier transform, there exists a constant $c$ such that $\int_\R \vert g(t)\vert^2\bigl(1-\cos(\alpha t)\bigr)\, dt \geq c\alpha^2$ for $\alpha$ sufficiently close to $0$.

\smallskip
\begin{proof}[Proof of Lemma \ref{bizarre}] For notational brevity, let us set $F:=\mathcal Ph$.

\smallskip
(a) We have $F(\alpha)=\langle h, \tau_\alpha h\rangle_{L_2(\R)}$, where $\tau_\alpha h(x):=h(x+\alpha)$, so $F$ is continuous by continuity of the map $\alpha\mapsto \tau_\alpha h$ from $\R$ into $L_2(\R)$. It is clear that $F$ is even and that $F(0)=\Vert h\Vert_2^2$. Finally, $F(\alpha)\leq \Vert h\Vert_2^2$ by the Cauchy-Schwarz inequality.

\smallskip
(b) Since $h\in W^{1,2}(\R)$, the map $\alpha\mapsto \tau_\alpha h$ is $\mathcal C^1\,$-$\,$smooth from $\R$ into $L^2(\R)$, with derivative 
$\alpha\mapsto \tau_\alpha h'$. Since $F(\alpha)=\langle h, \tau_\alpha h\rangle_{L_2(\R)}$, it follows that $F$ is $\mathcal C^1\,$-$\,$smooth on $\R$, with 
\begin{align*}
F'(\alpha)= \langle h, \tau_\alpha h'\rangle_{L_2(\R)}= \int_\R h(u)h'(u+\alpha)\, du=\int_\R h(u-\alpha)h'(u)\, du.
\end{align*}

Hence, the same argument shows that $F$ is in fact $\mathcal C^2\,$-$\,$smooth on $\R$, with 
\[ F''(\alpha)=-\int_\R h'(u-\alpha)h'(u)\, du.\]

We have $F'(0)=0$ because $F$ is even; 
and
\[ F''(0)=-\int_\R h'^2 <0.\]
The strict inequality holds because if we had $F''(0)=0$, then $h'$ would be equal to $0$ almost everywhere, so $h$ would be $0$  since it belongs to $W^{1,2}(\R)$.

\smallskip 
(c) Note that 
\[ F(\alpha)= \langle h, \tau_\alpha h\rangle_{L_2(\R)}=\frac12 \left( \Vert h\Vert_2^2+\Vert \tau_\alpha h\Vert^2-\Vert h-\tau_\alpha h\Vert_2^2\right)=\Vert h\Vert_2^2-\frac12 \Vert h-\tau_\alpha h\Vert_2^2,\]
so that $\Vert h\Vert_2^2-F(\alpha)=\frac12 \Vert h-\tau_\alpha h\Vert_2^2$. Now since $h\not\in W^{1,2}(\R)$, we have $\lim_{\alpha\to 0} \left\Vert \frac{\tau_\alpha h-h}{\alpha}\right\Vert=\infty$. Indeed, otherwise one could find a sequence $(\alpha_n)$ tending to $0$ such that $\frac{\tau_{\alpha_n} h-h}{\alpha_n}$ has a weak limit $u\in L_2(\R)$, and this would give that $h\in W^{1,2}(\R)$ with $h'=u$. This proves (c).

\smallskip (d) Since $F$ is an even function, it is enough to show that $F$ is right-differentiable at $0$ with $F'(0^+)<0$.

\smallskip Let $u_0<\dots <u_N$ be the discontinuity points of $h$. Also, let $I_0:=(-\infty, u_0)$, $I_{N+1}:=(u_n,\infty)$ and $I_k:=(u_{k-1}, u_{k})$ for $1\leq k\leq N$. Note that by assumption on $h$, the restriction of $h$ to each interval $I_k$ belongs to the Sobolev space $W^{1,2}(I_k)$. In particular, $h$ has a left limit $h(u_k^-)$ and a right limit $h(u_k^+)$ at each point $u_k$.

\smallskip
For every $\alpha>0$, we have
\[ F(\alpha)=\sum_{k=0}^N \int_{I_k} h(u) h(u+\alpha)\, du=: \sum_{k=0}^{N+1} F_k(\alpha).\]

We consider separately the functions $F_0$, $F_{N+1}$ and $F_k$ for $1\leq k\leq N$.

Let us start with $F_k$, $1\leq k\leq N$. For $0<\alpha<u_{k}-u_{k-1}$, we write 
\[ F_k(\alpha)=\int_{u_{k-1}}^{u_k-\alpha} h(u)h(u+\alpha)\, du + \int_{u_{k}-\alpha}^{u_k} h(u)h(u+\alpha)\, du=: F_{k,1}(\alpha)+F_{k,2}(\alpha).\]

Consider the open triangle $\Omega:=\{ (\alpha,\beta);\; u_{k-1}<\beta<u_k\;{\rm and}\; 0<\alpha<u_k-\beta\}$, and the function $G:\overline\Omega\to \R$ defined by $G(\alpha,\beta):= \int_{u_{k-1}}^\beta h(u)h(u+\alpha)\, du$. Since the restriction of $h$ to $I_k=(u_{k-1},u_k)$ belongs to the Sobolev space $W^{1,2}(I_k)$, the map $G$ is $\mathcal C^1$-$\,$smooth on $\Omega$, with $\partial_\alpha G(\alpha,\beta)=\int_{u_{k-1}}^\beta h(u)h'(u+\alpha)\, du$ and $\partial_\beta G(\alpha,\beta)= h(\beta)h(\beta+\alpha)$. Moreover, $G$ is continuous on $\overline\Omega$ and its partial derivatives extend continuously to $\overline{\Omega}$, with $\partial_\alpha G(\alpha, u_k-\alpha)=\int_{u_{k-1}}^{u_k-\alpha} h(u)h'(u+\alpha)\, du$ and $\partial_\beta G(\alpha, u_k-\alpha)=h(u_k-\alpha)h(u_k^-)$ for $0<\alpha<u_k-u_{k-1}$.  It follows that $F_{k,1}(\alpha)=G(\alpha, u_k-\alpha)$ is $\mathcal C^1$-$\,$smooth on $(0,u_{k}-u_{k-1})$, with $F_{k,1}'(\alpha)=\int_{u_{k-1}}^{u_k-\alpha} h(u)h'(u+\alpha)\, du - h(u_k-\alpha)h(u_k^-)$. To see this, observe that for any $\alpha_0\in (0,u_k-u_{k-1})$, the function $G_n(\alpha):= G(\alpha, u_k-\alpha -\frac1n)$ is well-defined in a neighbourhood of $\alpha_0$ if $n$ is large enough, that $G_n(\alpha)\to F_{k,1}(\alpha)$ uniformly in a neighbourhood of $\alpha_0$ as $n\to\infty$, and that $G'_n(\alpha)\to  \int_{u_{k-1}}^{u_k-\alpha} h(u)h'(u+\alpha)\, du - h(u_k-\alpha)h(u_k^-)$ uniformly in a neighbourhood of $\alpha_0$. Since $F_{k,1}'(\alpha)\to \int_{u_{k-1}}^{u_k} h(u)h'(u)\, du- h(u_k^-)^2$ as $0^+$, we deduce that $F_{k,1}$ is right-differentiable at $0$ with 
\begin{align*} F_{k,1}'(0^+)&= \int_{u_{k-1}}^{u_k} h(u)h'(u)\, du- h(u_k^-)^2\\
&=\left[\frac12 h(u)^2\right]_{u_{k-1}}^{u_k}- h(u_k^-)^2\\
&=-\frac12\bigl( h(u_k^-)^2+h(u_{k-1}^+)^2\bigr).
\end{align*}

Similar arguments show that $F_{k,2}(\alpha)=\int_{u_{k}-\alpha}^{u_k} h(u)h(u+\alpha)\, du$ is $\mathcal C^1$-$\,$smooth on $(0,c)$ for some $c>0$, with $F_{k,2}'(\alpha)=\int_{u_{k}-\alpha}^{u_k} h(u)h'(u+\alpha)\, du+h(u_k-\alpha)h(u_k^+)$; and it follows that $F_{k,2}$ is right-differentiable at $0$ with $F_{k,2}'(0^+)= h(u_k^-)h(u_k^+)$. Hence, $F_k$ is right-differentiable at $0$, with 
\[ F_k'(0^+)= h(u_k^-)h(u_k^+)-\frac12\,\bigl( h(u_k^-)^2+ h(u_{k-1}^+)^2\bigr).\]

One shows in the same way that $F_0(\alpha)=\int_{-\infty}^{u_0} h(u)h(u+\alpha)\, du$ and $F_{N+1}(\alpha)=\int_{u_N}^\infty h(u)h(u+\alpha)\, du$ are right-differentiable at $0$, with 
\[ F_0'(0^+)=h(u_0^-)h(u_0^+)-\frac12\, h(u_0^-)^2\qquad{\rm and}\qquad F_{N+1}'(0^+)=-\frac12\, h(u_N^+)^2.\]

Altogether, we see that $F=\sum_{k=0}^{N+1} F_k$ is right-differentiable at $0$, with 
\begin{align*}
F'(0^+)&=h(u_0^-)h(u_0^+)-\frac12\, h(u_0^-)^2\\
&\qquad +\sum_{k=1}^N \left( h(u_k^-)h(u_k^+)-\frac12\,\bigl( h(u_k^-)^2+ h(u_{k-1}^+)^2\bigr)\right) -\frac12\, h(u_N^+)^2\\
& =\sum_{k=0}^N h(u_k^-)h(u_k^+)-\frac12\sum_{k=0}^N \bigl( h(u_k^-)^2+ h(u_{k}^+)^2\bigr)\\
&=-\frac12\sum_{k=0}^N \left(h(u_k^+)-h(u_k^-) \right)^2<0.
\end{align*}
\epf


\smallskip We can now give the proof of Theorem \ref{gauss2}.
\begin{proof}[Proof of Theorem \ref{gauss2}] In what follows, we set $f:=\sqrt{\mathtt p}$ and $g:=\sqrt{\mathtt q}$. Also, for each $n\geq 1$, let 
\[ \lambda_n:=\frac{\uprod}{\vprod}\cdot\]

The following simple computation will be essential for the proof.
\begin{fact}\label{Psi} For every $n\geq 1$, we have 
\[ H(\mu_{\u,n},\mu_{\v,n})=\Psi(\lambda_n),\]
where $\Psi:(0,\infty)\to \R$ is the function defined by 
\[\Psi(\lambda):= \sqrt{\lambda} \int_\R f(\lambda t) g(t)\, dt.\]
\end{fact}
\begin{proof}[Proof of Fact \ref{Psi}] Denoting by $\tau$ the Lebesgue measure on $\R$, we have 
\[ \frac{d\mu_{\u,n}}{d\tau}(t)=\uprod\ \mathtt p(\uprod t)\qquad{\rm and} \qquad \frac{d\mu_{\v,n}}{d\tau}(t)=\vprod\ \mathtt q(\vprod t).\]

Hence
\begin{align*}
H(\mu_{\u,n},\mu_{\v,n})&=\int_\R\sqrt{\uprod}\,  f(\uprod t)\times \sqrt{\vprod} \, g(\vprod t)\, dt\\
&= \sqrt{\frac{\uprod}{\vprod}}\int_\R  f\left({\frac{\uprod}{\vprod}}\, s\right)\, g(s)\,  ds\\
&=\Psi (\lambda_n).
\end{align*}
\epf

\smallskip By Fact \ref{Psi} and Kakutani's Theorem, we know that $m_\u$ and $m_\v$ are not orthogonal if and only if
\begin{equation}\label{sumpsi} \sum_{n=1}^\infty (1-\Psi(\lambda_n))<\infty;
\end{equation}
and in that case, we have in particular that
\begin{equation}\label{limitpsi}
\Psi(\lambda_n)\to 1\quad\hbox{as $n\to\infty$}.
\end{equation}

\smallskip 
\begin{proof}[Proof of {\rm (1i)} in Theorem \ref{gauss2}] We first note that we must have $\varliminf \lambda_n>0$ and $\varlimsup \lambda_n<\infty$. Indeed, otherwise $B_\u$ and $B_\v$ are orthogonal with respect to measures not charging $\ker(e_0^*)$ by Proposition \ref{presque}. The measures 
$\mu_\u$ and $\mu_\v$ have this property since $\mu_{\u,0}$ and $\mu_{\v,0}$ do not charge $\{ 0\}$, so this cannot happen.

\smallskip Now, the key point is the following fact.

\begin{fact}\label{F1prime} The function $\Psi$ is continuous on $(0,\infty)$ with $0\leq \Psi(\lambda)\leq 1$, there is at most one $\lambda\in (0,\infty)$ such that $\Psi(\lambda)=1$, and if $\Psi(\lambda)=1$ then $\mathtt q(t)=\lambda \mathtt p(\lambda t)$ almost everywhere..
\end{fact}
\begin{proof}[Proof of Fact \ref{F1prime}]
 It is clear that $\Psi(\lambda)\geq 0$. Moreover, by the Cauchy-Schwarz inequality, we have  
  \[ \Psi(\lambda)\leq \left( \int_\R \lambda \mathtt p(\lambda t)\, dt\right)^{1/2}\left( \int_\R \mathtt q(t)\, dt\right)^{1/2}=1. \]
  
  If $\Psi(\lambda)=1$ then, by the equality case in Cauchy-Schwarz's inequality and since $\int_\R \mathtt p=\int_\R \mathtt q=1$, we must have $\mathtt q(t)=\lambda \mathtt p(\lambda t)$ almost everywhere. It follows easily that there can be at most one $\lambda$ such that $\Psi(\lambda)=1$. 
  Indeed, assume that $\Psi(\lambda)=1=\Psi(\lambda')$ for some $0<\lambda<\lambda'$. Then $\lambda \mathtt p(\lambda t)=\lambda' \mathtt p(\lambda' t)$ almost everywhere. Setting $c:=\lambda'/\lambda >1$, it follows that $\mathtt p(x)=c\, \mathtt p(c x)$ almost everywhere. Hence $\mathtt p(x)=c^n \mathtt p(c^n x)$ almost everywhere 
  for every  $n\in\Z$. So, with $I_n:=(c^n, c^{n+1}]$, we obtain that 
  \[ \int_{I_n} \mathtt p(x)\, dx=\int_{I_0} c^n \mathtt p(c^nx)\, dx=\int_{I_0} \mathtt p(x)\, dx\qquad\hbox{for all $n\in\Z$}.\]
  Since the intervals $I_n$ form a partition of $(0,\infty)$ and $\int_0^\infty \mathtt p(x)\, dx<\infty$, this implies that $\mathtt p(x)=0$ almost everywhere on $(0,\infty)$. Similarly $\mathtt p(x)=0$ almost everywhere on $(-\infty,0)$, hence $\mathtt p(x)=0$ almost everywhere on $\R$, which is a contradiction since $\int_\R \mathtt p=1$.
  
Finally, to prove that $\Psi$ is continuous, we note that $\Psi(\lambda)=\sqrt{\lambda} \, \langle f_\lambda, g\rangle_{L^2(\R)}$, where $f_\lambda(t):=f(\lambda t)$. Since $f\in L^2(\R)$, the map 
 $\lambda\mapsto f_\lambda$ is continuous from $(0,\infty)$ into $L^2(\R)$; 
 so $\Psi$ is indeed continuous. \epf 
 
 \smallskip It is now easy to conclude the proof of (1i).  Since $\Psi(\lambda_n)\to 1$ and $\Psi$ is continuous on $(0,\infty)$, any cluster point $\lambda\in (0,\infty)$ of the sequence $(\lambda_n)$ must be such that $\Psi(\lambda)=1$; hence, $(\lambda_n)$ has at most one cluster point in $(0,\infty)$ by Fact \ref{F1prime}. Since $0<\varliminf \lambda_n\leq\varlimsup\,\lambda_n<\infty$, it follows that $(\lambda_n)$ has a limit $a\in (0,\infty)$; and by Fact \ref{F1prime} again, we have  ($\Psi(a)=1$ and) $\mathtt q(t)=a \mathtt p(a t)$ almost everywhere.
  \epf

\smallskip
\begin{proof}[Proof of {\rm (1ii)} in Theorem \ref{gauss2}]  
 Let $\Phi:(0,\infty)\to \R$ and $\Theta:(0,\infty)\to \R$ be defined by 
\[ \Phi(\lambda):=\int_\R f(t)f(\lambda t)\, dt \qquad{\rm and}\qquad \Theta(\lambda):=\sqrt{\lambda}\, \Phi(\lambda).\]

Since $g(t)=\sqrt{a} f(a t)$ almost everywhere by (1i), we have 
\[ \Psi(\lambda)= \sqrt{\lambda} \int_\R f(\lambda t) \times \sqrt{a} f(a t)\, dt=\sqrt{\lambda/a}\, \Phi(\lambda/a)=\Theta(\lambda/a).\]

So, we see that $0\leq \Theta(\lambda)\leq 1$ for every $\lambda\in (0,\infty)$, and by (\ref{sumpsi}):
\begin{equation}\label{sumtheta}
\sum_{n=1}^\infty \bigl( 1-\Theta(\lambda_n/a)\bigr)<\infty.
\end{equation}

From that, we have to deduce that \begin{equation}\label{encoreune}\sum_{n=1}^\infty \left(1-\frac1a\, \lambda_n\right)^2<\infty.\end{equation}

To do this, we perform a change of variable in order to apply Lemma \ref{bizarre}. If we set $\lambda=e^\alpha$, then 
\begin{align*}
\Theta(\lambda)
&=e^{\alpha/2} \left( \int_0^{+\infty} f(t)f(e^\alpha t)\, dt+\int_{-\infty}^0\ f(t)f(e^\alpha t)\, dt\right)\\
&=e^{\alpha/2}\int_\R f(e^u)f(e^{\alpha+u})\, e^u du+ e^{\alpha/2}\int_\R f(-e^u)f(-e^{\alpha+u})\, e^u du\\
&=\int_\R h_+(u)h_+(u+\alpha)\, du+\int_\R h_-(u)h_-(u+\alpha)\, du,
\end{align*}
where $h_+$ and $h_-$ are the functions defined on $\R$ by 
\[ h_+(x):= f(e^x)e^{x/2}\qquad{\rm and}\qquad h_-(x):=f(-e^x)e^{x/2}.\]

Note that the functions $h_+$ and $h_-$ belong to $L_2(\R)$, and that 
\[ \Vert h_+\Vert_2^2+\Vert h_-\Vert_2^2=\int_\R f^2=1.\]

So, with the notation of Lemma \ref{bizarre}, we have 
\begin{equation}\label{applybizarre} \Theta(\lambda)=\mathcal Ph_+\bigl(\log(\lambda)\bigr)+\mathcal Ph_-\bigl(\log(\lambda)\bigr). \end{equation}

Since at least one of $h_\pm$ is non-zero, it follows that there exists a constant $c>0$ such that 
$1-\Theta(\lambda)\geq c\,\log(\lambda)^2$ {for $\lambda$ sufficiently close to $1$}. Hence $1-\Theta(\lambda)\geq c (\lambda-1)^2$ for $\lambda$ sufficiently close to $1$ (and some other constant $c>0$), so that (\ref{encoreune}) is indeed a direct consequence of (\ref{sumtheta}) since we already know that $\lambda_n\to a$. This concludes the proof of (1ii).
\epf

\begin{proof}[Proof of {\rm (2)} in Theorem \ref{gauss2}] By Kakutani's Theorem and with the notation of the proof of (1ii), we have to show that 
\[ \sum_{n=1}^\infty \bigl(1-\Theta(\lambda_n/a)\bigr)<\infty\quad\hbox{if and only if}\quad \sum\limits_{n=1}^\infty \left(1-\frac1a\, {\frac{\uprod}{\vprod}}\right)^2<\infty.\]

By (1ii), the ``only if'' implication is already known. To prove the converse, we keep the notation of the proof of (1ii). So we have as above
\[ \Theta(\lambda)=\mathcal Ph_+\bigl(\log(\lambda)\bigr)+\mathcal Ph_-\bigl(\log(\lambda)\bigr),\]
where the functions $h_+$ and $h_-$ are defined by 
\[ h_{\pm}(x)=f(\pm e^x)e^{x/2}.\]

Note that by assumption on $f$, the functions $h_+$ and $h_-$ are continuous on $\R$ and $\mathcal C^1\,$-$\,$smooth except at a finite number of points. Moreover, $h'_{\pm}\in L_2(\R)$. Indeed, we have $h_\pm'(x)=\pm f'(\pm e^x)e^{3x/2}+ \frac12 f(\pm e^x)e^{x/2}$ almost everywhere, and both terms belong to $L_2(\R)$ since
\[ \int_\R f'(\pm e^x)^2e^{3x} dx=\int_{\pm(0,\infty)} t^2f'(t)^2 dt<\infty\qquad{\rm and} \]
\[\int_\R f(\pm e^x)^2 e^x\, dx=\int_{\pm(0,\infty)} f(t)^2 dt <\infty.\]

So $h_{\pm}$ belongs to the Sobolev space $W^{1,2}(\R)$. Hence, by Lemma \ref{bizarre} and since $\Vert h_+\Vert_2^2+\Vert h_-^2\Vert_2=1$, we know that the function $F:=\mathcal Ph_++\mathcal Ph_-$ is $\mathcal C^2$-$\,$smooth on $\R$, with $F(0)=1$, $F'(0)=0$ and $F''(0)<0$. It follows that there is a constant $c>0$ such that $1-\Theta(\lambda)\sim c (\lambda-1)^2$ as $\lambda\to 1$; and this concludes the proof of the ``if'' implication in (2).
\epf

\begin{proof}[Proof of {\rm (3)} in Theorem \ref{gauss2}]  We keep the notation of the proofs of (1ii) and (2). This time, we have to show that 
\begin{equation}\label{sanslecarre}
\sum\limits_{n=1}^\infty \bigl(1-\Theta(\lambda_n/a)\bigr)<\infty\quad\hbox{if and only if}\quad \sum_{n=1}^\infty \left\vert 1-\frac1a\,\lambda_n\right\vert<\infty.
\end{equation}

As above, let $F:=\mathcal Ph_++\mathcal Ph_-$, so that $\Theta(\lambda)=F\bigl(\log(\lambda)\bigr)$. By assumption on 
$\mathtt p$, the functions $h_+(x)=f(e^x)e^{x/2}$ and $h_-(x)=f(-e^x)e^{x/2}$ are $\mathcal C^1$-$\,$smooth except at finitely many points, with $h'_{\pm}\in L_2(\R)$, and at least one of them has a discontinuity point. By Lemma \ref{bizarre}, it follows that there exists two constants $c>0$ and $c'<\infty$ such that $c\, \vert\alpha\vert \leq 1-F(\alpha)\leq c'\, \vert\alpha\vert$ for $\alpha$ sufficiently close to $0$. Hence $c\, \vert 1-\lambda\vert\leq 1-\Theta(\lambda)\leq c'\, \vert 1-\lambda\vert$ for $\lambda$ sufficiently close to $1$, and  (\ref{sanslecarre}) follows.
\epf

\smallskip
The proof of Theorem \ref{gauss2} is now complete. \end{proof}

\begin{remark} The proofs of parts (2) and (3) of Theorem \ref{gauss2} rely on the local properties of the function $\Theta(\lambda)=\sqrt{\lambda} \int_\R f(t)f(\lambda t)\, dt$ at 
$\lambda=1$. It is perhaps worth noticing that if $\Theta$ happens to be differentiable at $\lambda=1$, then we necessarily 
have $\Theta'(1)=0$. This is because $\Theta$ is ``symmetric with respect to $1$'', \textit{i.e.} $\Theta(1/\lambda)=\Theta(\lambda)$.
\end{remark}

\begin{remark} Part (3) of Theorem \ref{gauss2}  may be applied for example if $\mu_{\u,0}$ is the uniform distribution over some bounded interval $I\subseteq\R$  and $\mu_{\v,0}$ is the uniform distribution over $\frac1a\, I$.
\end{remark}

\begin{remark} Part (2) of Theorem \ref{gauss} remains valid under the following assumption on $f=\sqrt{p}$: the restrictions of $f$ to the intervals $(0,\infty)$ and $(-\infty,0)$ have locally integrable derivatives (in the distribution sense), and $tf'(t)\in L_2(\R)$.
\end{remark}
\bpf Looking back at the proof of (2), we see that the only thing to check is that the functions $h_+$ and $h_-$ belong to the Sobolev space $W^{1,2}(\R)$. We do this for $h_+$ and, for notational simplicity, we set $h:=h_+$ and we denote by $f$ the restriction of $f$ to $(0,\infty)$. So we have to show that if $f\in L_2(0,\infty)$ has a distributional derivative $f'\in L_1^{\rm loc}(0,\infty)$ such that $tf(t)\in L_2(0,\infty)$, then $h(x)=f(e^x)e^{x/2}\in W^{1,2}(\R)$. 

Since $f\in L_2(0,\infty)$, it is clear that $h\in L_2(\R)$. To find the distributional derivative of $h$, let us fix a test function $\varphi\in \mathcal D(\R)$. Let $\psi\in \mathcal D(0,\infty)$ be the function defined by the relation $\varphi(x)=\psi(e^x)e^{x/2}$, \textit{i.e.} $\psi(t)= \frac1{\sqrt{t}}\, \varphi(\log(t))$. Writing $\varphi'(x)= e^{3x/2} \psi'(e^x)+\frac12 e^{x/2}\psi(e^x)$ and using the change of variable $x=\log(t)$, an elementary  computation reveals that
\begin{align*}\int_\R h(x)\,\varphi'(x)\, dx&=-\int_0^\infty \left( e^{3x/2}f'(e^x)+\frac12\, e^{x/2}f(e^x)\right) \varphi(x)\, dx.
\end{align*}

This means that $h$ has a distributional derivative $h'\in L_1^{\rm loc}(\R)$ given (as expected) by $ h'(x)= e^{3x/2}f'(e^x)+\frac12 e^{x/2}f(e^x).$ Since $f\in L_2(0,\infty)$ and $tf'(t)\in L_2(0,\infty)$, we see that $h'\in L_2(\R)$, hence $h\in W^{1,2}(\R)$.

\epf


\subsection{Products of discrete measures}  Theorem \ref{gauss2} (3) shows that for weighted shifts acting on $\ell_p$, to admit non-orthogonal invariant measures which are products of absolutely continuous measures having some ``singularities'' is in fact a strictly stronger requirement on the weights than to admit equivalent 
invariant Gaussian product measures. The next proposition goes along the same lines, in an even more dramatic way. 
\bpr\label{gauss3} Let $\u$ and $\v$ be two weight sequences. 

\be
\item[\rm (a)] If $B_\u$ and $B_\v$ admit non-orthogonal invariant product measures $m_\u=\otimes_{n\geq 0} \mu_{\u,n}$ and $m_{\v}=\otimes_{n\geq 0}\mu_{\v,n}$ such that $\mu_{\u,0}$ and $\mu_{\v,0}$ have non-zero discrete parts and do not charge $\{ 0\}$, then $\left\vert\frac{\uprod}{\vprod}\right\vert$ is eventually constant {\rm (}and hence $B_\u$ and $B_\v$ admit equivalent invariant Gaussian product measures{\rm )}. 
\item[\rm (b)] Assume that $B_\u$ and $B_\v$ admit equivalent non-trivial invariant product measures  or, more generally, non-trivial invariant product measures $m_\u=\otimes_{n\geq 0} \mu_{\u,n}$ and $m_{\v}=\otimes_{n\geq 0}\mu_{\v,n}$ such that $\mu_{\u,n}$ and $\mu_{\v,n}$ have the same support for each $n\geq 0$.  
Let $S:={\rm supp}(\mu_{\u,0})={\rm supp}(\mu_{\v,0})$. If either $S$ is compact or $0\not\in S$, then 
$\left\vert\frac{\uprod}{\vprod}\right\vert=1$ for all $n\geq 1$; and hence $\u=\v$ if $\u$ and $\v$ are positive weight sequences.
\ee
\epr
\bpf For each $n\geq 1$, let us set as usual $\lambda_n:=\frac{\uprod}{\vprod}\cdot$

\smallskip
(a) Since the measures $\mu_{\u,n}$ and $\mu_{\v,n}$ do not charge $\{ 0\}$, one can find a countable multiplicative subgroup $\mathcal S$  of $\K\setminus\{ 0\}$ such that the discrete parts of all measures $\mu_{\u,n}$ and $\mu_{\v,n}$ are supported on $\mathcal S$ and $\mathcal S$ contains the set $\{ u_n;\; n\geq 1\}\cup\{ v_n;\; n\geq 1\}$. Denote by $\tau_{\mathcal S}$ the counting measure on $\mathcal S$. Let also $\tau_c$ be a continuous measure such that the continuous parts of all measures $\mu_{\u,n}$ and $\mu_{\v,n}$ are absolutely continuous with respect to $\tau_c$. Then, one can take $\tau:= \tau_c+\tau_{\mathcal S}$ to compute 
$H(\mu_{\u,n},\mu_{\v,n})$, for every $n\geq 0$.

For each $n\geq 0$, write 
\[ \mu_{\u,n}=\alpha_n \tau_c+\mathtt p_n \tau_{\mathcal S}\qquad {\rm and} \qquad \mu_{\v,n}=\beta_n \tau_c+\mathtt q_n \tau_{\mathcal S}\]
for some non-negative measurable functions $\alpha_n, \mathtt p_n, \beta_n, \mathtt q_n$. Without loss of generality, we may assume that $\alpha_n \mathbf 1_{\mathcal S}= \beta_n \mathbf 1_{\mathcal S}=0$ and $\mathtt p_n \mathbf 1_{\K\setminus\mathcal S}= \mathtt q_n \mathbf 1_{\K\setminus\mathcal S}=0$, so that we also have 
\[ \mu_{\u,n}=(\alpha_n+\mathtt p_n) \tau\qquad{\rm and} \qquad \mu_{\v,n}=(\beta_n+\mathtt q_n) \tau.\]
Hence,
\[ H(\mu_{\u,n},\mu_{\v,n})=\int_\K \sqrt{\alpha_n+\mathtt p_n}\,\sqrt{\beta_n+\mathtt q_n}\, d\tau = \langle\, \sqrt{\alpha_n+\mathtt p_n}, \sqrt{\beta_n+\mathtt q_n}\; \rangle_{L_2(\tau)}.\]

Now, by Kakutani's Theorem we know that $H(\mu_{\u,n},\mu_{\v,n})\to 1$ as $n\to\infty$. Since $\Vert \sqrt{\alpha_n+\mathtt p_n}\Vert_{L_2(\tau)}= \Vert \sqrt{\beta_n+\mathtt q_n}\Vert_{L_2(\tau)}=1$, it follows that \[ \bigl\Vert \sqrt{\alpha_n+\mathtt p_n}-\sqrt{\beta_n+\mathtt q_n}\,\bigr\Vert_{L_2(\tau)}\to0.\] Moreover, by assumption on $\alpha_n,\beta_n, \mathtt p_n, \mathtt q_n$, we have 
\[ \bigl\Vert \sqrt{\alpha_n+\mathtt p_n}-\sqrt{\beta_n+\mathtt q_n}\bigr\Vert_{L_2(\tau)}^2=\bigl\Vert \sqrt{\alpha_n}-\sqrt{\beta_n}\bigr\Vert_{L_2(\tau_c)}^2+\bigl\Vert \sqrt{\mathtt p_n}-\sqrt{\mathtt q_n}\bigr\Vert_{L_2(\tau_{\mathcal S})}^2.\]

Hence, we have in particular that $\Vert \sqrt{\mathtt p_n}-\sqrt{\mathtt q_n}\Vert_{L_2(\tau_{\mathcal S})}^2\to 0$ as $n\to\infty$; and since $\tau_{\mathcal S}$ is the counting measure on $\mathcal S$, this means that 
\begin{equation}\label{l2} \sum_{s\in\mathcal S} \Bigl(\, \sqrt{\mathtt p_n(s)}-\sqrt{\mathtt q_n(s)}\,\Bigr)^2\to 0\qquad\hbox{as $n\to\infty$}.\end{equation}

So far, we have not used the fact that the measures $m_\u$ and $m_\v$ are invariant under $B_\u$ and $B_\v$. Set $\mathtt p:=\mathtt p_{0}$ and $\mathtt q:= \mathtt q_0$. Then, the discrete parts of $\mu_{\u,0}$ and $\mu_{\v,0}$ are respectively $\sum_{s\in\mathcal S} \mathtt p(s)\, \delta_s$ and $\sum_{s\in\mathcal S} \mathtt q(s)\, \delta_s$. So, by the invariance properties and since $\mathcal S$ is a multiplicative group containing $\uprod$ and $\vprod$, the discrete parts of $\mu_{\u,n}$ and $\mu_{\v,n}$ (for $n\geq 1$) are respectively $\sum_{s\in\mathcal S} \mathtt p(\uprod s)\, \delta_s$ and  $\sum_{s\in\mathcal S} \mathtt q(\vprod s)\, \delta_s$; in other words, we have $\mathtt p_n(s)= \mathtt p(\uprod s)$ and $\mathtt q_n(s)=\mathtt q(\vprod s)$ for all $s\in\mathcal S$. Hence, if we set $f:=\sqrt{\mathtt p}$ and $g:=\sqrt{\mathtt q}$, (\ref{l2}) can be re-written as follows: 
\[ \sum_{s\in\mathcal S} \bigl( f(\uprod s)-g(\vprod s)\bigr)^2\to 0\qquad\hbox{as $n\to\infty$}\]
or, equivalently,
\[ \sum_{s\in\mathcal S} \bigl( f(\lambda_n s)-g(s)\bigr)^2\to 0.\]

It follows in particular, $f(\lambda_ns)\to g(s)$ for all $s\in\mathcal S$. So, taking any $s$ such that $g(s)> 0$ (such an $s$ exists since the discrete part of $\mu_{\v,0}$ is non-zero), we see that one can find $\varepsilon >0$, e.g. $\varepsilon:= g(s)/2$, such that $ f(\lambda_n s)\geq \varepsilon$ for all $n$ sufficiently large. It follows that the sequence $(\lambda_n)$ can take only finitely many distinct values: indeed, otherwise $ f(s')\geq \varepsilon$ for infinitely many $s'\in\mathcal S$, a contradiction since $f\in\ell_2(\mathcal S)$. Now, assume that $\vert \lambda_n\vert$ is not eventually constant. Then, since $\lambda_n$ takes only finitely many values, $\lambda_n\in\mathcal S$ and $f(\lambda_n s)\to g(s)$ for all $s\in\mathcal S$, we see that one can find $\lambda,\lambda'\in\mathcal S$ with $\vert \lambda\vert \neq\vert\lambda'\vert$ such that $f(\lambda s)=f(\lambda's)$ for all $s\in\mathcal S$. Setting $\alpha:=\lambda/\lambda'$, we then have $f(\alpha^k s)=f(s)$ for any $s\in\mathcal S$ and all $k\in\N$; and taking any $s$ such that $f(s)> 0$ (again, such an $s$ exists), we obtain a contradiction since $f\in\ell_2(\mathcal S)$ and the $\alpha^k$, $k\in\N$ are pairwise distinct. This concludes the proof of (a).

\smallskip
(b)  Let us fix $n\geq 1$ and, towards a contradiction, assume that $\vert\lambda_n\vert\neq 1$, say $\vert\lambda_n\vert<1$.

By the invariance properties, we have ${\rm supp}(\mu_{\u,n})=\frac1\uprod \, S$ and ${\rm supp}(\mu_{\v,n})=\frac1\vprod \, S$, so that $\lambda_n S=S$. Hence $\lambda_n^r S=S$ for all $r\in\Z$. If $S$ is compact, it follows that $S=\{ 0\}$ by letting $r\to +\infty$;  hence $\mu_{\u,0}=\delta_0=\mu_{\v,0}$, which is the required contradiction. If $0\not\in S$, then we get $S=\emptyset$ by letting $r\to -\infty$, which is another contradiction.
\epf

\bco Let $\u$ and $\v$ be positive weight sequences such that $\sum_{n=1}^\infty \frac{1}{(\uprod)^p}<\infty$ and $\sum_{n=1}^\infty \frac{1}{(\vprod)^p}<\infty$. The weighted shifts $B_\u$ and $B_\v$ acting on $\ell_p$ are not orthogonal with respect to product measures whose marginals have non-zero discrete parts and do not charge $\{ 0\}$ if and only if $\frac\uprod\vprod$ is eventually constant; and in that case they admit equivalent invariant product measures with purely discrete marginals not charging $\{ 0\}$.
\eco
\bpf The ``only if'' part follows from Proposition \ref{gauss3} (a). Conversely, assume that $\lambda_n:=\frac\uprod\vprod$ is eventually constant, say $\lambda_n=\lambda$ for $n>n_0$. Let $\mathcal S$ be the multiplicative subgroup of $(0,\infty)$ generated by the set 
$\{ u_n;\; n\geq 1\}\cup\{ v_n;\; \geq 1\}$. Let $\mathtt p:\mathcal S\to \R_+$ be a strictly positive probability density function such that $\sum_{s\in\mathcal S} \mathtt p(s)\, s^p<\infty$, and let $\mathtt q:\mathcal S\to \R_+$ be the probability density function defined by $\mathtt q(s):= \mathtt p(\lambda s)$. Finally, let 
$\mu_{\u, 0}:=\sum_{s\in\mathcal S} \mathtt p(s)\delta_s$ and $\mu_{\v,0}:=\sum_{s\in\mathcal S} \mathtt q(s)\delta_s$, and denote by $m_\u=\oplus_{n\geq 0} \mu_{\u,n}$ and $m_\v=\oplus_{n\geq 0}\mu_{\v,n}$ the associated $B_\u\,$-$\,$invariant and $B_\v\,$-$\,$invariant measures on $\K^{\Z_+}$. By Corollary \ref{charge} and the assumption on $\u$, $\v$ and $\mathtt p$, the measures $m_\u$ and $m_\v$ are supported on $\ell_p$, and they have purely discrete marginals not charging $\{ 0\}$. Finally, by the proof of Proposition \ref{gauss3} (a), we have $H(\mu_{\u,n},\mu_{\v,n})=\sum_{s\in\mathcal S} \sqrt{\mathtt p(\lambda_n s)\mathtt q(s)}$ for all $n\geq 0$, so that 
$H(\mu_{\u,n},\mu_{\v,n})=\sum_{s\in\mathcal S} \mathtt q(s)=1$ for all $n>n_0$. Since $\mu_{\u,n}\sim\mu_{\v,n}$ for all $n\geq 0$, it follows that $m_\u\sim m_\v$. This concludes the proof.
\epf




\begin{remark} In Proposition \ref{gauss3} (b), one cannot replace ``equivalent'' by ``non-orthogonal''. For example, let $\u$ and $\v$ be defined as follows: $u_n=2$ for all $n\geq 1$, $v_1=1$, $v_2=4$ and $v_n=2$ for all $n\geq 3$. Then $\u\neq \v$ (!), and yet $B_\u$ and $B_\v$ admit non-orthogonal product measure $m_\u=\otimes_{n\geq 0}\mu_{\u,n}$ and $m_\v=\otimes_{n\geq 0}\mu_{\v,n}$ for which $\mu_{\u,0}=\mu_{\v,0}$ has a compact support not containing $0$; for example, $\mu_{\u,0}=\mu_{\v,0}$ could be the uniform distribution  on the interval $[1,3]$.
\end{remark}

\section{Additional facts}\label{add}
\subsection{Invariant measures with symmetry properties} 
It is clear that if an operator $T$ acting on a Polish topological vector space $X$ admits a non-trivial invariant measure $m$, then the measure $\widetilde m$ defined by $\widetilde m(A):= \frac12(m(A)+m(-A))$ is a non-trivial invariant measure for both $T$ and $-T$, and that $\widetilde m$ is additionally \emph{symmetric}, \textit{i.e.} 
$\widetilde m(-A)=\widetilde m(A)$ for every Borel set $A\subseteq X$. The next proposition goes along the same lines.

\bpr\label{symm} Let $X$ be a Polish space, and let $T$ be a continuous self-map of $X$ admitting an invariant measure $m$. Let also $\mathbf G$ be a compact abelian group acting continuously on $X$, and assume that $Tg=\beta(g) T$ for all $g\in \mathbf G$, where $\beta:\mathbf G\to \mathbf G$ is a continuous map invariant under the Haar measure of $\mathbf G$. Let $\widetilde m$ be the measure on $X$ defined by 
\[ \widetilde m(A):=\int_{\mathbf G} m(g^{-1}A)\, dg,\]
where $dg$ denotes integration with respect to the Haar measure of $\mathbf G$. Then $\widetilde m$ is $\mathbf G\,$-$\,$invariant and both $gT\,$-$\,$invariant and $Tg\,$-$\,$invariant for every $g\in \mathbf G$. 
\epr
\bpf It is clear that $\widetilde m$ is $\mathbf G\,$-$\,$invariant. Moreover, if $h\in \mathbf G$ then, for any bounded Borel function $f:X\to \R^+$, we have 
\begin{align*}
\int_X f\circ(\beta(h)T)\, d\widetilde m&=\int_X\left( \int_{\mathbf G} f(\beta(h)T( gx))\, dg\right)dm(x)\\
&=\int_X \left( \int_{\mathbf G} f(\beta(h)\beta(g) T(x))\, dg\right) dm(x)\\
&=\int_X \left( \int_{\mathbf G} f(\beta(h)g T(x))\, dg\right) dm(x)\\
&=\int_X \left( \int_{\mathbf G} f(g T(x))\, dg\right) dm(x)\\
&=\int_{\mathbf G} \left(\int_X f(g u)\, dm(u)\right) dg \\
&=\int_X f\, d\widetilde m.\\
\end{align*}
 Since $\beta$ is necessarily onto (because the Haar measure of $\mathbf G$ has full support), it follows that $\widetilde m$ if $gT\,$-$\,$invariant for every $g\in \mathbf G$; and hence also $Tg\,$-$\,$invariant since $Tg=\beta(g)T$.
  \epf

\bco Let $X$ be a Polish space, and let $T_1,T_2$ be two continuous self-maps of $X$. Let also $\mathbf G$ be a compact abelian group acting continuously on $X$, and assume that $T_ig=gT_i$ for all $g\in\mathbf G$ and $i=1,2$. Finally, let $x_0\in X$ be a fixed point for the action of $\mathbf G$.  If $T_1$ and $T_2$ admit equivalent {\rm (}resp. non-orthogonal{\rm )} invariant measures $m_1, m_2$ not charging $\{ x_0\}$, then they also admit equivalent {\rm (}resp. non-orthogonal{\rm )} invariant measures not charging $\{ x_0\}$ which are additionally $\mathbf G\,$-$\,$invariant.
\eco
\bpf Let $\widetilde{m_1}$ and $\widetilde{m_2}$ be defined as in Proposition \ref{symm}. It is clear that $\widetilde{m_1}$, $\widetilde{m_2}$ do not charge $\{ x_0\}$. So we just have to check that if $m_1\sim m_2$ then $\widetilde{m_1}\sim\widetilde{m_2}$, and that if $\widetilde{m_1}\perp\widetilde{m_2}$ then $m_1\perp m_2$.

Assume that $m_1\sim m_2$. Let $A\subseteq X$ be a Borel set such that $\widetilde{m_1}(A)=0$. Then $m_1(gA)=0$ for almost every $g\in \mathbf G$; hence $m_2(gA)=0$ for almost every $g\in \mathbf G$ since $m_2\ll m_1$, and hence $\widetilde{m_2}(A)=0$. This shows that $\widetilde{m_2}\ll \widetilde{m_1}$; and similarly $\widetilde{m_1}\ll \widetilde{m_2}$.

Assume that $\widetilde{m_2}\perp \widetilde{m_1}$. Let $A\subseteq X$ be a Borel set such that $\widetilde{m_1}(A)=0$ and $\widetilde{m_2}(X\setminus A)=0$. Then $m_1(gA)=0$ for almost every $g\in\mathbf G$ and $m_2(X\setminus gA)=m_2(g(X\setminus A))=0$ for almost every $g\in\mathbf G$. So one can find at least one $g$ such that $m_1(gA)=0=m_2(X\setminus gA)$, which shows that $m_1\perp m_2$.
\epf

\bco Let $X$ be a Polish topological vector space, and let $T\in\mathcal L(X)$. Assume that $\kappa T$ admits a non-trivial invariant measure for some $\kappa\in\K$. Then, there exists a non-trivial measure which is $aT\,$-$\,$invariant for all $a\in\K$ such that $\vert a\vert=\vert\kappa\vert$.
\eco
\bpf Apply Proposition \ref{symm} with the group $\mathbf G:=\{ \omega \in \K;\; \vert \omega\vert=1\}$ acting by multiplication on $X$.
\epf

\bco\label{unitary} Let $B_\u$ be a backward shift acting on $\ell_p$. If $B_\u$ admits a non-trivial invariant measure, then there exists a measure $m$ on $\ell_p$ with full support which is invariant for all backward shifts $B_{\v}$ such that $\vert v_n\vert=\vert u_n\vert$ for all $n\geq 1$.
\eco
\bpf Since $\sum_{n=1}^\infty\frac1{\vert\uprod\vert^p}<\infty$, we know that $B_\u$ admit an invariant measure $m$ with full support (which can even be taken to be Gaussian and ergodic for $B_\u$). Let \[ \mathbf G:=\left\{ (\omega_n)_{n\geq 0}\in \K^{\Z_+};\; \vert \omega_n\vert=1\;\;\hbox{for all $n\geq 0$}\right\}.\] This is a compact abelian group acting continuously on $X:=\K^{\Z_+}$ by coordinatewise multiplication. Moreover, if we set $T:= B_\u$ acting on $X$, we have $B_\u g=\sigma(g) B_\u$ for every $g\in\mathbf G$, where $\sigma:\mathbf G\to \mathbf G$ is (the restriction to $\mathbf G$ of) the canonical backward shift. Hence, considering $m$ as a Borel measure on $X$, we may apply  Proposition \ref{symm}. This gives a measure $\widetilde m$ on 
$\K^{\Z_+}$ which is $B_\v\,$-$\,$invariant for every weight sequence $\v$ such that $\vert v_n\vert=\vert u_n\vert$ for all $n\geq 1$. Moreover, since $B_\u$ is supported on $\ell_p$ and $\ell_p$ is $\mathbf G\,$-$\,$invariant, $\widetilde m$ is supported on $\ell_p$; and since $m$ has full support, it is readily checked that $\widetilde m$ has full support. 
\epf

\begin{remark} It is well-known (see \cite{Sh}) that two weight sequences $\u$ and $\v$ are such that $\vert u_n\vert=\vert v_n\vert$ for all $n\geq 1$ if and only if the backward shifts $B_\u$ and $B_\v$ (acting on any $\ell_p$) are unitarily similar, \textit{i.e.} there exists an isometry $J$ of $\ell_p$ such that $B_\v=J B_\u J^{-1}$.
\end{remark}

\subsection{Product measures charging small subspaces} 
The following proposition says in essence that whether or not a product measure $\K^{\Z_+}$ invariant under a weighted shift $B_\w$ is supported on some ``small'' subspace of $\K^{\Z_+}$ depends on the rate of growth of the products $\wprod$.
\bpr Let $\mu_0$ be a Borel probability measure on $\K$, let $\w$ be a  sequence of non-zero scalars, and let $\mu_\w=\otimes_{n\geq 0}\mu_n$ be the $B_\w$-$\,$invariant measure on $\K^{\Z_+}$ defined by $\mu_0$, \textit{i.e.} $\mu_n(A)=\mu_0(\wprod A)$ for each $n\geq 1$ 
and every Borel set $A\subseteq\K$. 

\be
\item[\rm (a)] Let $p\in [1,\infty)$. If $\int_\K \vert t\vert^p d\mu_0(t)<\infty$ and $\sum_{n=1}^\infty \frac1{\vert \wprod\vert^p}<\infty$, then $\mu_\w(\ell_p)=1$.
\item[\rm (b)] If there exists a summable sequence of positive real numbers $(\varepsilon_n)$ such that  $$\sum_{n=1}^\infty \mu_0 \bigl(\vert t\vert >\vert\wprod\vert\,   \varepsilon_n\bigr)<\infty,$$ then $\mu_\w(\ell_1)=1$. 
\item[\rm (c)] If $\vert\wprod\vert\not\to\infty$ and $\mu_0\neq \delta_0$, then $\mu_\w(c_0)=0$. If $\vert\wprod\vert\not\to\infty$ and $\mu_0$ does not have compact support, then $\mu_\w(\ell_\infty)=0$. 
\ee
\epr
\bpf (a) This is Corollary \ref{charge}.

\smallskip
(b) Since $\ell_1$ is a tail subset of $\K^{\Z_+}$, we have $\mu_\w(\ell_1)=0$ or $1$, by Kolmogorov's $0\,$-$\,1$ law. Let 
\[ A:=\bigl\{ x\in \K^{\Z_+};\; \vert x_n\vert\leq \varepsilon_n\quad\hbox{for all $n\geq 1$}\bigr\}.\]

Obviously $A\subseteq\ell_1$, and we have 
\begin{align*} \mu_\w(A)&=\prod_{n=1}^\infty \mu_n\bigl(\vert t\vert \leq \varepsilon_n\bigr)\\
&=\prod_{n=1}^\infty \mu_0\bigl(\vert t\vert \leq \vert\wprod\vert\, \varepsilon_n\bigr).
\end{align*}

Since $\sum_{n=1}^\infty \bigl(1- \mu_0(\vert t\vert \leq \vert\wprod\vert\, \varepsilon_n)\bigr)<\infty$ by assumption, it follows that $\mu_\w(A)>0$, and hence $\mu_\w(\ell_1)=1$.

\smallskip
(c) Assume that $\mu_\w(c_0)>0$. Then, there exists a compact set $K\subseteq c_0$ such that $\mu_\w( K) >0$. By the well-known description of the compact subsets of $c_0$, one can find a sequence of positive real numbers $(\varepsilon_n)$ tending to $0$ such that 
\[ K\subseteq \bigl\{ x\in \K^{\Z_+};\; \vert x_n\vert \leq \varepsilon_n\quad\hbox{for all $n\geq 0$}\bigr\}.\]

As in the proof of (b), it follows that 
\[ \sum_{n=1}^\infty \mu_0\bigl(\vert t\vert > \vert\wprod\vert\, \varepsilon_n\bigr)<\infty.\]

In particular, $\mu_0(\vert t\vert > \vert\wprod\vert\, \varepsilon_n)\to 0$ as $n\to\infty$. If $\mu_0\neq \delta_0$, this implies that no subsequence of $(\vert\wprod\vert \varepsilon_n)$ can tend to $0$; and hence $\vert\wprod\vert\to\infty$.

\smallskip
The second part of (c) is proved in a similar way.
\epf

\bco For any probability measure $\mu_0$ on $\K$, one can find a sequence of positive real numbers $\w$ such that $\mu_\w(\ell_1)=1$.
\eco
\bpf This is clear by (b): choose a sequence of positive numbers $(X_n)$ such that $\sum_{n=1}^\infty \mu_0(\vert t\vert >X_n)<\infty$, and then take $\w$ such that $2^{-n}\wprod\geq X_n$ for all $n\geq 1$.
\epf

\bco If  the measure $\mu_0$ is such that $\sum_{n=1}^\infty \mu_0(\vert t\vert >C^n)<\infty$ for some constant $C$, then  $\mu_\w(\ell_1)=1$ for any weight sequence $\w$ such that $\varliminf \vert \wprod\vert^{1/n} >C$. In particular,  if $\mu_0(\vert t\vert >X)=O\bigl(\log(X)^{-\alpha}\bigr)$ as $X\to\infty$ for some constant $\alpha>1$, then $\mu_\w(\ell_1)=1$ for any weight sequence $\w$ such that $\varliminf \vert w_n\vert^{1/n} >1$.
\eco
\bpf Take $\varepsilon_n:= (C/\rho)^n$ in (b), where $\rho:= \varliminf \vert w_n\vert^{1/n} $.
\epf

\subsection{When frequently hypercyclic vectors are the same} As mentioned in Remark \ref{same}, the following result is due to S. Charpentier and the third author. We thank S. Charpentier for allowing us to include it here.
\bpr
Let $\u$ and $\v$ be two weight sequences such that $B_\u$ and $B_\v$ are frequently hypercyclic on $\ell_p$. If $\frac{\uprod}{\vprod}$ has a non-zero limit as $n\to\infty$, then $B_\u$ and $B_\v$ have the same frequently hypercyclic vectors.
\epr
\begin{proof}
In what follows, we set $a:=\lim\limits_{n\to\infty} \frac{\uprod}{\vprod}\in \K\setminus\{ 0\}.$

\smallskip
Let $x$ be a frequently vector for $B_{\u}$; we want  to show that $x$ is also a frequenty hypercyclic vector for $B_\v$. So, we fix  $y=\sum_{k=0}^{d}y_k e_k\in c_{00}$ and $\varepsilon>0$, and our task is to show that the set $\mathcal N_{B_\v}\bigl( x, B(y,\varepsilon)\bigr)$ has positive lower density. 

\smallskip Since $x$ is a frequently hypercyclic vector for $B_\u$, it is enough to find a vector $z\in \ell_p$ and $\alpha >0$ such that any sufficiently large $n\in \mathcal N_{B_\u}\bigl( x, B(z, \alpha)\bigr)$ belongs to $\mathcal N_{B_\v}\bigl( x, B(y,\varepsilon)\bigr)$. 
We consider 
\[z:=\sum_{k=0}^{d} a\, \frac{v_1\cdots v_k}{u_1\cdots u_k}\, y_k\, e_k\in c_{00}\]
and 
\[ \alpha:= \varepsilon/2K\qquad{\rm where}\quad K:=\sup_{n,k\ge 0}\left|\frac{v_{k+1}\cdots v_{k+n}}{u_{k+1}\cdots u_{k+n}}\right|.\] 
 
 Note that $K<\infty$  since for every $n,k\ge 0$, we have
\[\left|\frac{v_{k+1}\cdots v_{k+n}}{u_{k+1}\cdots u_{k+n}}\right|=\left|\frac{u_{1}\cdots u_{k}}{v_{1}\cdots v_{k}}\right| \left|\frac{v_{1}\cdots v_{k+n}}{u_{1}\cdots u_{k+n}}\right|\le \left(\sup_{m\geq 1} \left\vert\frac{u_1\cdots u_m}{v_1\cdots v_m}\right\vert\right)\left(\sup_{m\geq 1}\left\vert \frac{v_1\cdots v_m}{u_1\cdots u_m}\right\vert\right).\]

Let us show that if $n\in \mathcal N_{B_{\u}}\bigl(x,B(z,\frac{\varepsilon}{2K})\bigr)$ is large enough, then $n\in\mathcal N_{B_\v}\bigl( x, B(y,\varepsilon)\bigr)$. We have

\begin{align*}
\|B_{\v}^nx-y\|^p&=\sum_{k=0}^{d}|(v_{k+1}\cdots v_{k+n})x_{k+n}-y_k|^p+\sum_{k=d+1}^{\infty} |(v_{k+1}\cdots v_{k+n})x_{k+n}|^p\\
&\le \sum_{k=0}^{d}\left|\frac{v_{k+1}\cdots v_{k+n}}{u_{k+1}\cdots u_{k+n}}\right|^p\left|(u_{k+1}\cdots u_{k+n})x_{k+n}-\frac{u_{k+1}\cdots u_{k+n}}{v_{k+1}\cdots v_{k+n}}\, y_k\right|^p\\
&\qquad\qquad+\sum_{k=d+1}^{\infty} \left|\frac{v_{k+1}\cdots v_{k+n}}{u_{k+1}\cdots u_{k+n}}\right|^p |(u_{k+1}\cdots u_{k+n})x_{k+n}|^p\\
&\le K^p\|B_{\u}^nx-\tilde{z}^{(n)}\|^p\\
\end{align*}
where
\[\tilde{z}^{(n)}:=\sum_{k=0}^{d} \frac{u_{k+1}\cdots u_{k+n}}{v_{k+1}\cdots v_{k+n}}y_k\,e_k.\]
 
 Therefore, if $n\in \mathcal N_{B_{\u}}\bigl(x,B(z,\frac{\varepsilon}{2K})\bigr)$, we get
\[
\|B_{\v}^nx-y\|\le K\|B_{\u}^nx-z\|+K\|z-\tilde{z}^{(n)}\|\le \frac{\varepsilon}{2}+K\|z-\tilde{z}^{(n)}\|;
 \]
and the desired result follows because
\begin{align*}
\|z-\tilde{z}^{(n)}\|&\le \sum_{k\le d} \left|z_k-\frac{u_{k+1}\cdots u_{k+n}}{v_{k+1}\cdots v_{k+n}}\,y_k\right|\\
&\le \|y\|_{1}\max_{k\le d}\left| a\, \frac{v_1\cdots v_k}{u_1\cdots u_k}-\frac{u_{k+1}\cdots u_{k+n}}{v_{k+1}\cdots v_{k+n}}\right|\\
&\le \|y\|_{1}\left(\max_{k\le d}\left|\frac{v_1\cdots v_k}{u_1\cdots u_k}\right|\right)\left(\sup_{m\ge n}\left|a -\frac{u_{1}\cdots u_{m}}{v_{1}\cdots v_{m}}\right|\right)\xrightarrow[n\to \infty]{} 0.
\end{align*}
\end{proof}

\subsection{Using Shepp's Theorem} In this section, we show how results like Shepp's theorem from \cite{Shepp} mentioned before Theorem \ref{gauss2} can be used in the context of weighted shifts. Specifically, we will make use of the following theorem. For any measure $m$ on $X=\R^{\Z_+}$ and 
$\pmb\alpha\in X$, let us denote by $m_{\pmb\alpha}$ the translate of $m$ by $\pmb\alpha$, which is the measure on $X$ defined by \[ m_{\pmb\alpha}(A):= m(A+\pmb\alpha).\] 
  
 \bth\label{known} Let $\widetilde\mu$ and $\widetilde \mu'$ be two measures on $\R$, and let $\widetilde m$ and $\widetilde m'$ be the product measures on $\R^{\Z_+}$ with marginals 
 $\widetilde\mu_n:=\widetilde \mu$ and $\widetilde\mu_n':=\widetilde \mu'$.  Let also $\pmb\alpha=(\alpha_n)_{n\geq 0}\in \R^{\Z_+}$, and assume that the measures $\widetilde m$ and 
 $\widetilde m'_{\pmb\alpha}$ are not orthogonal.
 \be
 \item[\rm (1)] If $\widetilde\mu=\widetilde\mu'$, then $\sum_{n=0}^\infty \alpha_n^2<\infty$.
 \item[\rm (2)] If $\widetilde\mu$ and $\widetilde\mu'$ have a moment of order $2$, then $\sum_{n=0}^\infty (\alpha_n-\alpha)^2<\infty$ for some $\alpha\in\R$.
 \ee
 \eth
 \bpf (1) is (the first third of) \cite{Shepp}*{Theorem 1}.
 
 \smallskip To prove (2), we use a method devised by Dudley \cite{Dudley} in order to generalize Shepp's theorem.  Denote by $e_n^*$, $n\geq 0$ the coordinate linear functionals on $\Omega:=\R^{\Z_+}$. The assumption on $\widetilde\mu$ and $\widetilde\mu'$ implies that $e_n^*$ belongs to $L_2(\Omega, \widetilde m)\cap L_2(\Omega,\widetilde m')$ with $L_2\,$-$\,$norms 
 respectively equal to $\int_\R x^2\, d\widetilde\mu(x)$ and $\int_\R x^2\, d\widetilde\mu'(x)$, hence independent of $n$. Likewise, $\int_\Omega e_n^* \, d\widetilde m=\int_\R x\, d\widetilde\mu(x)$ and $\int_\Omega e_n^* \, d\widetilde m'=\int_\R x\, d\widetilde\mu'(x)$ do not depend on $n$.
 
 \smallskip Assume first that  $\int_\Omega e_n^* \, d\widetilde m=0$ and $\int_\Omega e_n^* \, d\widetilde m'=0$ for all $n\geq 0$. We claim that in this case, we have $\sum_{n=0}^\infty \alpha_n^2<\infty$.   Towards a contradiction, assume that $\sum_{n=0}^\infty \alpha_n^2=\infty$. Then, one can find a sequence of real numbers $(\beta_n)$ such that $\sum_{n=0}^\infty \beta_n^2<\infty$, $\beta_n\alpha_n\geq 0$ for all $n$ and $\sum_{n=0}^\infty\beta_n\alpha_n=\infty$. Now, consider the measure $\tau:=\frac12(\widetilde m+\widetilde m')$. By our assumption, the sequence $(e_n^*)_{n\geq 0}$ is orthogonal and bounded in $L_2(\Omega,\tau)$. Since $\sum_{n=0}^\infty \beta_n^2<\infty$, it follows that the series $\sum \beta_n e_n^*$ is convergent in $L_2(\Omega,\tau)$. Hence, there is an increasing sequence of integers $(N_k)_{k\geq 0}$ such that if we define $f_k:=\sum_{n=0}^{N_k} e_n^*$, then the sequence of linear functionals $(f_k)$ converges $\tau\,$-$\,$a.e. on $\Omega$. So, the linear subspace
 \[ E:=\bigl\{ x\in \Omega;\; \hbox{the sequence $(f_k(x))$ is convergent}\bigr\}\subseteq \Omega\] 
 is such that $\tau(E)=1$ and hence $\widetilde m(E)=1$. On the other hand, we have $f_k(\pmb\alpha)=\sum_{n=0}^{N_k} \beta_n\alpha_n$, so the sequence $(f_k(\pmb\alpha))$ is \emph{not} convergent, \textit{i.e.} $\pmb\alpha\not\in E$; and since $E$ is a linear subspace of $\Omega$, it follows that 
 $E\cap (E+\pmb\alpha)=\emptyset$. So we have $\tau(E+\pmb\alpha)=0$, and hence $\widetilde m_{\pmb\alpha}'(E)=\widetilde m'(E+\pmb\alpha)=0$. Since $\widetilde m$ and $\widetilde m'_{\pmb\alpha}$ are not orthogonal, this is a contradiction.
 
 \smallskip Now, let us consider the general case. Let $c:=\int_\R x\, d\widetilde\mu(x)$ and $c':=\int_\R x\, d\widetilde\mu'(x)$, so that $\int_\Omega e_n^*\, d\widetilde m=c$ and $\int_\Omega e_n^*\, d\widetilde m'=c'$ for all $n\geq 0$. Then $\int_\Omega e_n^* \, d\widetilde m_{\mathbf c}=0$ and $\int_\Omega e_n^* \, d\widetilde m'_{\mathbf c'}=0$ for all $n\geq 0$, where $\mathbf c:=(c,c,\dots)$ and $\mathbf c':=(c',c',\dots)$. Since the measures $\widetilde m_{\mathbf c}$ and $\bigr(\widetilde m'_{\mathbf c'}\bigl)_{\pmb\alpha+\mathbf c-\mathbf c'}=\bigl(\widetilde m'_{\pmb\alpha}\bigr)_{\mathbf c}$ are not orthogonal, it follows that $\sum_{n=0}^\infty (\alpha_n-(c'-c))^2<\infty$. 
 \epf

\smallskip From Theorem \ref{known}, it is essentially a formal matter to deduce the following theorem, which is not far from saying that if two weighted shifts are not orthogonal with respect to product measures whose marginals do not charge $\{ 0\}$, then they admit equivalent invariant Gaussian product measures.
\bth\label{Shepp} Let $\u$ and $\v$ be two sequences of non-zero scalars, and let $m_\u=\otimes_{n\geq 0} \mu_{n,\u}$ and $m_\v=\otimes_{n\geq 0} \mu_{n,\v}$ be two product measures on $\K^{\Z_+}$, respectively $B_\u\,$-$\,$invariant and $B_\v\,$-$\,$invariant, with $\mu_{0,\u}(\{ 0\})=0=\mu_{0,\v}(\{ 0\})$. Assume that the measures 
$m_\u$ and $m_\v$ are not orthogonal. Finally, let $\lambda_n:=\frac{\uprod}{\vprod}$ for all $n\geq 1$.  

\smallskip\be
\item[\rm (1)] If there exists $a\in\K\setminus\{ 0\}$ such that $\mu_{0,\v}(A)=\mu_{0,\u}(a A)$ for every Borel set $A\subseteq \K$, then $\sum_{n=1}^\infty \bigl(1- \vert a^{-1}\lambda_n\vert\bigr)^2<\infty$.
\item[\rm (2)] If $\int_{\K\setminus\{ 0\}} \bigl( \log \vert t\vert\bigr)^2 d\mu_{0,\u}(t)<\infty$ and $\int_{\K\setminus\{ 0\}} \bigl( \log \vert t\vert\bigr)^2 d\mu_{0,\v}(t)<\infty$, then there exists $a\in\K\setminus\{ 0\}$ such that $\sum_{n=1}^\infty \bigl(1- \vert a^{-1}\lambda_n\vert\bigr)^2<\infty$.
\ee
\eth
\bpf (1) Since $\mu_{0,\u}(\{ 0\})=0= \mu_{0,\v}(\{ 0\})$, we view $\mu_{0,\u}$ and $\mu_{0,\v}$ as measures on $G:=\K\setminus\{ 0\}$, and hence we view $m_\u$ and $m_\v$ as measures on $G^{\Z_+}$. Note that $G^{\Z_+}$ is an abelian group under entry-wise multiplication.

\smallskip Consider the measures $m=\otimes_{n\geq 0} \mu_n$ and $m'=\otimes_{n\geq 0} \mu'_{n}$ where $\mu_n:=\mu_{0,\u}$ and $\mu'_n:=\mu_{0,\v}$ for all $n\geq 0$. 
If we set $\mathbf U:=(1,u_1, u_1u_2,\dots)$, $\mathbf V:=(1, v_1, v_1v_2,\dots )$ and $\mathbf a:=(a,a,a,\dots)$, then we have for every Borel set $B\subseteq G^{\Z_+}$:
\[ m_\u(B)=m(\mathbf UB)\qquad{\rm and}\qquad m_\v(B)=m'(\mathbf VB)=m(\mathbf a\mathbf V B).\]

\smallskip
Now, consider the map $L :G^{\Z_+}\to \R^{\Z_+}$ defined by $L(t_0, t_1,\dots ):=\bigl( \log\vert t_0\vert, \log\vert t_1\vert,\dots \bigr)$. Denote by $\widetilde m_\u$, $\widetilde m_\v$, $\widetilde m$ and $\widetilde m'$ the images of $m_\u$, $m_\v$, $m$
 and $m'$ under this map $L$. Then $\widetilde m=\otimes_{n\geq 0} \widetilde\mu_n$ and $\widetilde m'=\otimes_{n\geq 0} \widetilde\mu'_n$, where $\widetilde \mu_n$ and 
 $\widetilde\mu'_n$ are the images of $\mu_{0,\u}$ and $\mu_{0,\v}$ under the map $t\mapsto\log\vert t\vert$, for every $n\geq 0$. Moreover, since $L$ is a group homomorphism, we have for every Borel set $B\subseteq \R^{\Z_+}$:
 \[ \widetilde m_\u(B)=\widetilde m\bigl(B+L(\mathbf U)\bigl)\qquad{\rm and}\qquad \widetilde m_\u(B)=\widetilde m\bigl( B+L(\mathbf{aV})\bigl).\]


 
 \smallskip Now, the measures $\widetilde m_\u=m_\u\circ L^{-1}$ and $\widetilde m_\v=m_\v\circ L^{-1}$ are not orthogonal since $m_\u$ and $m_\v$ are not orthogonal. So their translates by $-L(\mathbf{aV})$ are not orthogonal either.  In other words, if we set
 \[ \pmb\alpha:= L(\mathbf{U}\mathbf{V}^{-1}\mathbf{a}^{-1})=(0, \log\vert a^{-1}\lambda_1\vert, \log\vert a^{-1}\lambda_2\vert,\dots ),\]
 then the measure $\widetilde m$ and $\widetilde m_{\pmb\alpha}$ are not orthogonal. By Theorem \ref{known} (1), it follows that $\sum_{n=0}^\infty \alpha_n^2<\infty$, \textit{i.e.}
 \[ \sum_{n=1}^\infty (\log\vert a^{-1}\lambda_n\vert)^2<\infty;\]
 which concludes the proof.
 
 \smallskip
 (2) The proof is similar, using Theorem \ref{known} (2).
\epf

\bco\label{bien2} Let $B_\u$ and $B_\v$ be two weighted shifts acting on $\ell_p$. The following are equivalent.
\be
\item[\rm (a)] $B_\u$ and $B_\v$ are not orthogonal with respect to product measures whose marginals are absolutely continuous with respect to Lebesgue measure on $\K$.
\item[\rm (b)] $B_\u$ and $B_\v$ are not orthogonal with respect to product measures whose marginals do not charge $\{ 0\}$ and are such that $\log\vert t\vert\in L_2$.
\item[\rm (c)] There exists $\kappa>0$ such that $\sum_{n=1}^\infty \left( 1-\kappa\left\vert \frac{\uprod}{\vprod}\right\vert \right)^2<\infty$.
\item[\rm (d)] $B_\u$ and $B_\v$ admit equivalent invariant Gaussian product measures.
\ee 
\eco
\bpf By Theorem \ref{Shepp} and keeping the same notation, we just have to check that if $\mu_{0,\u}$ and $\mu_{0,\v}$ are absolutely continuous with respect to Lebesgue measure  on $\K$, then there exists $a\in \K\setminus\{ 0\}$ such that $\mu_{0,\v}(A)=\mu_{0,\u}(aA)$ for every Borel set $A\subseteq\K$. When $\K=\R$ this is done in the proof of (1i) in Theorem \ref{gauss2}. We repeat the argument here, in a way that makes it work when $\K=\C$ as well. Write $\mu_{0,\u}=\mathtt p(t)dt$ and $\mu_{0,\v}=\mathtt q(t)dt$. Then $\Psi(\lambda_n)\to 1$ as $n\to\infty$, where $\Psi:\K\setminus\{ 0\}\to [0,1]$ is the function defined by $\Psi(\lambda):=\vert \lambda\vert^{d/2}\, \int_\K f(\lambda t) g(t)\, dt$ and $d=1$ or $2$ depending on whether $\K=\R$ or $\C$. Moreover, since $B_\u$ and $B_\v$ must be similar by Proposition \ref{presque}, the sequence $(\lambda_n)$ has a cluster point $a\in \K\setminus\{ 0\}$. Then $\Psi(a)=1$ by continuity of $\Psi$, and it follows that $\mathtt q(t)=\vert a\vert^{d} \,\mathtt p(at)$ almost everywhere.\epf


\subsection{Kakutani's Theorem}  In this section, our aim is to prove the following variant of Kakutani's Theorem from \cite{Kak}.

\bth\label{kakutbis} Let $\nu=\otimes_{n\geq 0}\nu_n$ and $\nu'=\otimes_{n\geq 0} \nu'_n$ be two product probability measures on $\Omega=\prod_{n\geq 0} \Omega_n$. The measures $\nu$ and $\nu'$ are orthogonal if and only if $\prod_{n=0}^\infty H(\nu_n,\nu'_n)=0$. Moreover, if $\prod_{n=0}^\infty H(\nu_n,\nu'_n)>0$ and $\nu_n\ll\nu'_n$ for all $n$, then $\nu\ll\nu'$.
\eth

\smallskip  Note that this variant was used in the proofs of Theorem \ref{gauss2} and Proposition \ref{gauss3}. We give a rather detailed proof since we were not able to locate one in the literature. However, this is essentially a ``copy and paste'' of Kakutani's original proof of his theorem.

\smallskip

\begin{proof}[Proof of Theorem \ref{kakutbis}]
For each $n\geq 0$, let us choose a Radon-Nikodym derivative of $\nu_n$ with respect to $\nu_n'$, \textit{i.e.} a measurable function $\mathtt p_n:\Omega_n\to \R^+$ such that 
\[ \nu_n=\mathtt p_n\, \nu'_n+ \alpha_n,\]
where $\alpha_n$ is a positive measure orthogonal to $\nu_n'$. Let also $f_n:= \sqrt{\mathtt p_n}$. With this notation, we have 
\[ H(\nu_n,\nu'_n)=\int_{\Omega_n} \sqrt{\mathtt p_n}\, d\nu'_n=\int_{\Omega_n} f_n\, d\nu'_n.\]

\medskip 
(i) Assume that $\prod_{n=0}^\infty H(\nu_n,\nu'_n)>0$, and let us show that the measures $\nu$ and $\nu'$ are not orthogonal, with $\nu\ll\nu'$ if $\nu_n\ll\nu'_n$ for all $n$.

For each $k\geq 0$, let $F_k:=\otimes_{n=0}^k f_n$, \textit{i.e.} $F_k:\Omega\to \R^+$ is the function defined by 
\[ F_k(\omega):=\prod_{n=0}^k f_n(\omega_n).\]

We claim that $F_k$ converges in $L_2(\Omega, \nu')$ to some $F:\Omega\to \R^+$. Indeed, if $0\leq p<q$, then
\begin{align*}
\Vert F_q-F_p\Vert_{L_2(\nu')}^2&= \Vert F_p\Vert^2_{L_2(\nu')}\left\Vert 1-\otimes_{n=p+1}^q f_n\right\Vert^2_{L_2(\nu')}\\
&\leq \left\Vert 1-\otimes_{n=p+1}^q f_n\right\Vert^2_{L_2(\nu')}\\
&\leq 2\left( 1- \prod_{n=p+1}^q H(\nu_n,\nu'_n)\right);
\end{align*}
so we see that $(F_k)_{k\geq 0}$ is a Cauchy sequence in $L_2(\Omega,\nu)$.

By $L_2\,$-$\,$convergence and since $\int_{\Omega_n} f_n^2\, d\nu'_n\geq \left(\int_{\Omega_n} f_n\, d\nu'_n\right)^2$ for all $n\geq 0$, we have
\[ \int_\Omega F^2\, d\nu'=\lim_{k\to\infty} \int_\Omega F_k^2\, d\nu'\geq \prod_{n=0}^\infty H(\nu_n,\nu'_n)^2>0.\]

Moreover, for any cylinder set $A=[A_0,\dots ,A_N]\subseteq \Omega$, we have 
\begin{align*}
\nu(A)&=\prod_{n=0}^N \nu_n(A_n)\\
&\geq \prod_{n=0}^N \int_{A_n} \mathtt p_n\, d\nu'_n\\
&= \prod_{n=0}^N \int_{A_n}  f_n^2\, d\nu'_n \geq \int_A F^2\, d\nu',
\end{align*}
where the last inequality follows from the fact that $\int_{\Omega_n} f_n^2 \, d\nu'_n\leq 1$ for every $n> N$. 

Therefore, the positive measure $F^2\, \nu'$ is non-zero and absolutely continuous with respect to both $\nu'$ and $\nu$, and hence $\nu$ and $\nu'$ are not orthogonal. Finally, if $\nu_n\ll \nu'_n$ for all $n$, \textit{i.e.} $\alpha_n=0$, then we see that $\nu(A)=\int_A F^2\, d\nu'$ for all cylinder sets $A$, so $\nu\ll \nu'$.

\smallskip
(ii) Now, assume that $\prod_{n=0}^\infty H(\nu_n,\nu'_n)=0$, and let us show that $\nu\perp \nu'$. It is enough to show that for any $\varepsilon >0$, one can find a measurable set $E\subseteq \Omega$ such that $\nu'(E)<\varepsilon$ and $\nu(\Omega\setminus E)<\varepsilon$.

Let $\varepsilon >0$, and let us choose $N\geq 0$ such that $\prod_{n=0}^N H(\nu_n,\nu'_n)<\varepsilon $, \textit{i.e.}
\[ \int_\Omega F_N\, d\nu'<\varepsilon.\]
 Let us also choose, for each $n\in\llbracket 0,N\rrbracket$, a measurable set $A_n\subseteq \Omega_n$ such that $\nu'_n(A_n)=0$ and $\alpha_n(\Omega_n\setminus A_n)=0$, and let 
 \[ E:=\{ \omega\in\Omega;\; F_N(\omega)\geq 1\}\cup \{ \omega\in\Omega;\; \omega_n\in A_n\;\hbox{for some $n\in \llbracket 0,N\rrbracket$}\}.\]
  We have on the one hand $\nu'(E)=\nu'(F_N\geq 1)<\varepsilon$; and, on the other hand,
 \begin{align*}
 \nu(\Omega\setminus E)&=\otimes_{n=0}^\infty (\mathtt p_n\, \nu'_n+\alpha_n)\, \bigl( \{ F_N<1\}\cap [\Omega_0\setminus A_0,\dots ,\Omega_N\setminus A_N]\bigr)\\
 &= \otimes_{n=0}^N (\mathtt p_n\, \nu'_n)\bigl(\{ (\omega_0,\dots ,\omega_N);\; \sqrt{\mathtt p_0(\omega_0)\cdots \mathtt p_N(\omega_N)}<1\}\bigr)\\
 &= \int_{\sqrt{\mathtt p_0\otimes\cdots \otimes \mathtt p_N}<1} \mathtt p_0\otimes\cdots \otimes \mathtt p_N\, d\nu'\\
 &\leq\int_\Omega \sqrt{\mathtt p_0\otimes\cdots \otimes \mathtt p_N}\, d\nu'= \int_{\Omega} F_N\, d\nu' <\varepsilon.
 \end{align*}
 
This concludes the proof.
\epf

\section{Some questions} Of course, the main open question that remains at the end of this work is the one from which we started.
\begin{question} Find a characterization of the pairs of weight sequences $(\u,\v)$ such that the weighted shifts $B_\u$ and $B_\v$ acting on $\ell_p$ are orthogonal. 
\end{question}

\smallskip In Section \ref{Sec3}, we saw that orthogonality might  be related to non-similarity, in a way which is not yet clear. In particular, in view of Theorem \ref{presquebis}, the following question is natural.
\begin{question}\label{q2} Is it true that if $B_\u$ and $B_\v$ admit non-trivial invariant measures and are orthogonal, then either 
\[\varliminf_{n\to\infty}\;\max_{0\leq d\leq N} \,\left\vert\frac{u_1\cdots u_{n+d}}{v_1\cdots v_{n+d}}\right\vert=0\qquad\hbox{for all $N\geq 0$} \]
or
\[ \varlimsup_{n\to\infty} \;\min_{0\leq d\leq N} \,\left\vert\frac{u_1\cdots u_{n+d}}{v_1\cdots v_{n+d}}\right\vert=\infty \qquad\hbox{for all $N\geq 0$}?\]
\end{question}

\smallskip
In the same spirit, Example \ref{non} shows that non-similarity does not imply orthogonality. Conversely, it would be quite interesting to know whether similarity implies non-orthogonality.
\begin{question}\label{q3} Suppose that $B_\u$ and $B_\v$ admit non-trivial invariant measures and are similar, \textit{i.e.}
\[ 0<\varliminf \,\left\vert \frac{\uprod}{\vprod}\right\vert\leq\varlimsup \,\left\vert \frac{\uprod}{\vprod}\right\vert<\infty.\]

Does it follow that $B_\u$ and $B_\w$ are non-orthogonal?
\end{question}

\smallskip Note that  a positive answer to Question \ref{q2} would imply a positive answer to Question \ref{q3}. A very special case of the latter has a positive answer by Corollary \ref{unitary}: if $B_\u$ and $B_\v$ admit non-trivial invariant measures and are \emph{unitarily} similar, then they are non-orthogonal.

\smallskip
Next, in view of Theorem \ref{gauss2}, Proposition \ref{gauss3} and Theorem \ref{Shepp}, it is natural to ask whether the existence of equivalent non-trivial invariant product measures always implies the existence of  equivalent invariant \emph{Gaussian} product measures.
\begin{question} Is it true in general that if $B_\u$ and $B_\v$ admit equivalent non-trivial invariant product measures $m_\u=\otimes_{n\geq 0} \mu_{\u,n}$ and $m_\v=\otimes \mu_{\v,n}$, then they admit equivalent invariant Gaussian product measures? 
\end{question}

\smallskip From a more general point of view, the following question seems also natural.
\begin{question} If $T_1$ and $T_2$ are non-orthogonal operators, does it follow that they admit equivalent non-trivial invariant measures? 
\end{question}

\smallskip In the same spirit, one may ask
\begin{question} When do two weighted shifts $B_\u$ and $B_\v$ share a common non-trivial invariant measure?
\end{question}

\smallskip Finally, the following seems to be unknown.
\begin{question} Let $B_\w$ be a weighted shift acting on $c_0(\Z_+)$. Is it true that $B_\w$ admits a non-trivial invariant measure if and only if $\vert\wprod\vert\to\infty$ as $n\to\infty$?
\end{question}





\end{document}